\documentclass[10pt] {article} 
  \usepackage{amsmath}
    \usepackage{amssymb}
  \usepackage[dvipsnames]{xcolor}   
     \usepackage{pdfsync} 
 
\usepackage[title]{appendix}

\newtheorem{theorem}{Theorem}[section]
\newtheorem{lemma}{Lemma}[section]

\newcommand{\eqnsection}{
   \renewcommand{\theequation}{\thesection.\arabic{equation}}
   \makeatletter
   \csname @addtoreset\endcsname{equation}{section} 
   \makeatother}
   
\def \ov{\overline}

\def \be{\begin{equation}}
\def \ee{\end{equation}}
\def \bt{\begin{theorem}} 
\def \et{\end{theorem}}
\def \bl{\begin{lemma}} 
\def \el{\end{lemma}}
\def \bea{\begin{eqnarray}}
\def \eea{\end{eqnarray}}
\def \bas{\begin{eqnarray*}}
\def \eas{\end{eqnarray*}}

\def \al{\alpha}
\def \bb{\beta}
\def \ga{\gamma}
\def\Ga{\Gamma}
\def \de{\delta}
\def \De{\Delta}
\def \ep{\epsilon}

\def \la{\lambda}

\def \om{\omega}
\def \Om{\Omega}

\def \si{\sigma}

\def \th{\theta}

\def \ze{\zeta}

\def \ff{\infty}
\def \wh{\widehat}
\def \wt{\widetilde}

\def \CC{{\cal C}}

\def \FF{{\cal F}}

\def\b1{\mathbf 1}
\def \({\left(}
\def \){\right)}

\def \nn{\nonumber}
 
\def \Proof{\noindent{\bf Proof $\,$ }}

\def \bc{\begin{center} }
\def \ec{\end{center} }
\def \bs{\begin{slide} }
\def \es{\end{slide} }

\def\square{{\vcenter{\vbox{\hrule height.3pt
        \hbox{\vrule width.3pt height5pt \kern5pt
           \vrule width.3pt}
        \hrule height.3pt}}}}
\def\qed{{\hfill $\square$ \bigskip}}

\eqnsection

\begin{document}

\title{Exact moduli of continuity for the local times of Feller Brownian motions }

 \author{   P.J. Fitzsimmons\,\,\, Jay Rosen   }
\maketitle

  \begin{abstract}

We examine the  modulus of continuity, in the spatial variable, of the local time process of Feller Brownian motion (FBM) on the half-line $[0,\infty)$.
Briefly, a FBM is a strong Markov process on $[0,\infty)$ that moves like standard Brownian motion on $[0,\infty)$ up until it first encounters the state $0$. The process returns to $(0,\infty)$, either continuously (like reflecting Brownian motion) or by jumping to a (random) positive state chosen according to a specified measure.   The present work is a continuation and application of our earlier work \cite{FMR} with Michael Marcus on the moduli of continuity for the local times of a Markov process built by piecing together (``rebirthing") the paths of another Markov process with finite lifetime.
We first establish a general result on the resolvent and local times for a rebirthed process with a special holding state (the state $0$ for FBM). We show how our  earlier approach \cite{FMR} using the Eisenbaum Isomorphism Theorem on an assemblage of excursions  works out in this context.  This knowledge is then used as an approximation device to obtain our main result on the exact uniform moduli of continuity for the local time of FBM on a spatial interval of the form $(0,1]$. Extensions are made, under certain conditions, to the more delicate situation of the spatial interval $[0,1]$. We also consider briefly the case of more general diffusions on $[0,\infty)$. 

\end{abstract}

\section{Introduction}

In 1963, K. Ito and H.P. McKean \cite{IM1} examined in detail the infinitesimal generators and paths of a class  of strong Markov processes on $[0,\infty)$ whose  generators had been described in 1952 by W. Feller \cite{F}. Briefly, a process is in this class if, when started at $x>0$,   it behaves like a standard Brownian motion until it reaches the state $0$. Returns to $(0,\infty)$ from $0$ can occur continuously or via jumps, and the boundary point $0$ may be ``sticky'',  the process accumulating a positive amount of time there. Perhaps the best way to understand these processes is by using excursion theory, and this story is told in detail in \cite{R83}, \cite{R84}, and \cite{Blu}.

In this paper we confine our attention to the recurrent case.  A recurrent Borel right process $(\wt X_{t}, \wt P^{x})$ with state space $[0,\ff)$ is called a Feller Brownian motion if it behaves like ordinary Brownian motion up until the hitting time $\si$ of the point $0$. Such a process admits local times, and our goal in this paper is to describe the modulus of continuity of the local time, in the spatial variable. We draw on results and ideas from our earlier work \cite{FMR} concerning such matters for rebirthed Markov processes. 
The novelty of the present paper lies in dealing with the fact that a Feller Brownian motion  may have  infinitely many excursions from $0$ --- including those that start with a jump from $0$ to a state $x>0$ --- in each finite time interval.

A given Feller Brownian motion $\wt X$ has a $1$-potential operator $U^{1}$ which, restricted to $(0,\ff)$, has a continuous density $u ^1 (x,u)$ with respect to  Lebesgue measure. Let $\wt L^{u}_{s}$ denote the local times of $\wt X$, normalized so that for
 $y, u>0$,
 \begin{equation}
  \wt E^{ y}\(\int_{0}^{\ff}  e^{-  s}\, d_{s}\wt L^{u}_{s}\)=   u ^1 (y,u). \label{int.00}
\end{equation}

 Our  main goal in this paper is to prove the following Theorem.
 
 \bt\label{theo-uniformmodulus} For any $y>0$
   \be
\lim_{h\to 0}\sup_{\stackrel{|u-v|\le h }{ u,v\in (0,1] }} \frac{|\wt L_{t}^{u}-\wt L_{t}^{v}|}{  (|u-v| \log  1/| u-v|)^{1/2} }= \sup_{u\in (0,1]}\(2\wt L_{t}^{u}\)^{1/2},  \quad a.e.\,\,\, t,  \,\,\, \wt P^{y}\,\,
a.s.\label{int.1}
\ee
 \et

In Section 2 we establish, for later use, a general resolvent decomposition for rebirthed Markov processes. We also revisit the work in \cite{FMR} concerning an exact uniform modulus of continuity for the local times of a rebirthed process.  The present situation is complicated by the need to allow the special state $0$ to be a holding point. (This need arises when approximating a Feller Brownian motion by a simpler process, in Section 3.) 

In Section 3 we state and prove the main result,   Theorem \ref{theo-uniformmodulus}, concerning the modulus of continuity for the local times of a Feller Brownian motion.

By \cite[Theorem 3]{R83} with $p_{1}=0$, $p_{2}=q/\sqrt{2}$, $p_{3}=p$ and $p_{4}=\eta$, we have 
\begin{equation}
U^{\la}h(x)=V^{\la}h(x) +e^{-\sqrt{2\la}\,\, x}\,\,\frac{ph(0)+ q\sqrt{2}\int_{0}^{\ff} e^{-\sqrt{2\la}\,\, y}h(y)\,dy +\eta V^{\la}h}{p\la+q\sqrt{\la}+\int_{0}^{\ff}\(1-e^{-\sqrt{2\la}\,\, x}\)\,d\eta(x)},\label{r.42r}
\end{equation}
where $ V^{\la}$ is the resolvent for Brownian motion killed the first time it hits $0$,
\[\sqrt{2}\int_{0}^{\ff} e^{-\sqrt{2\la}\,\, y}h(y)\,dy=\sqrt{\la}\, R^{\la}h(0)\] where $ R^{\la}$ is the resolvent for reflected Brownian motion and \[\int_{(0,\ff)}\(1-e^{-x}\)\,d\eta(x)<\ff.\]
$p$ measures the stickiness of $0$, $q$ measures the tendency for $\wt X$ to exit $0$
continuously, and $\eta $ is the L\'evy measure of jumps of $\wt X$ from $0$ into $(0,\ff)$. 

In Section 5 we consider the problem of extending Theorem \ref{theo-uniformmodulus} to $u,v\in [0,1]$. When  $p\neq 0$   the $1$-potential operator $U^{1}$ does not have  a continuous density with respect to  Lebesgue measure on  $[0,\ff)$ . However, when  $p=0$ but $q>0$,  $U^{1}$  does  have  a continuous density $u ^1 (y,u)$  with respect to  Lebesgue measure on $[0,\ff)$ and  in Section 5
we present a large class of examples where Theorem \ref{theo-uniformmodulus} can be extended to $u,v\in [0,1]$.

When $p=0$ and  $q=0$ we have
$u ^1 (y,0)=0$, so we have to normalize the local times of $\wt X$  differently. In this case we have not been able to find exact uniform moduli of continuity for the local times, but we provide examples where the local times, suitably normalized, are  jointly continuous on $[0,1]\times R_{+}^{1}$.

 In Section 4 we indicate how  the main result, Theorem \ref{theo-uniformmodulus}, can be extended to the more general context in which the Brownian motion on $(0,\infty)$ is replaced by a more general regular diffusion.

 \section{ Local times of rebirthed Markov processes with holding parameter $\mathbf{J}$ and jumping in measure $\mathbf{\mu}$}\label{sec-Borel}

In Theorem \ref{theo-mborel} we obtain a general result about the potential densities of   fully rebirthed Borel right processes: 
    Let $S$ a be locally compact space with a countable base.  
Let   $X=(\Om,  \FF_{t}, X_t,\th_{t}, P^x)$ be a transient Borel right process with state space $S$  and  continuous symmetric strictly positive  $p-$potential   densities  
\be
v^{\la}=\{v^{\la}(x,y), x,y\in S \}, \qquad \la\geq 0,
\ee
 with respect to some $\si$--finite  positive measure $m$ on $S$. Let $\ze=\inf\{t\,|\,X_t=\De \}$,  where $\De$ is the cemetery state for   $X$  and assume       that $\ze<\ff$   almost surely.   
As usual, we use $E^x$ for expectations with respect to $P^x$.

    We use the following notation.   For measurable functions  $h$ on $S$ we  define  the $\la$--th potential operator,
\begin{equation} \label{6.1mm}
V^{\la}h(x)=\int_Sv^{\la}(x,y)h(y)\,dm(y). 
\end{equation}  
As a probabilistic statement this is 
\begin{equation} \label{6.2mm}
V^{\la}h(x)=  E^{x}\( \int_{0}^{\ff}e^{-\la t}h\(X_{t}\) \,dt \)=  E^{x}\( \int_{0}^{\ze}e^{-\la t}h\( X_{t}\) \,dt \),
\end{equation}
and in particular
\begin{equation} \label{6.2mm1}
V^{\la}1(x)=  E^{x}\( \int_{0}^{\ze}e^{-\la t}  \,dt \)=\frac{1}{\la} E^{x}\(1-e^{-\la \ze}   \).
\end{equation}

 For positive measures $\mu$  on $S$ we define,  
  \begin{equation} 
  \mu h =\int_S   h(x)\,d\mu(x). 
\end{equation}

  Note that by (\ref{6.1mm}),
\begin{equation} \label{6.5mm}
  \mu V^\la h=\int_S \int_S v^{\la }(x,y)h(y)\,dm(y)\,d\mu(x).
\end{equation}
In particular,     when $h\equiv 1 $, and using (\ref{6.2mm1}), (\ref{6.5mm}) is
\begin{equation} \label{6.5mmma}
  \mu V^\la 1=\int_S \int_S v^{\la }(x,y)\,dm(y)\,d\mu(x)=\frac{1}{\la} E^{\mu}\(1-e^{-\la \ze}\).
\end{equation}

  Let $\mu$ be a probability measure on $S$ and let $0$ be a point not in $S$. We       modify $X$ so that   at the end of its lifetime it immediately goes to $0$ and stays there an independent exponential holding time $J$ with parameter $\bb$, after which it is ``reborn'' with probability measure $\mu $.  (I.e.,  $\wt  X_J\in B\subset S$ with probability $\mu(B) $, after which  it continues to evolve the way $X$ did, being reborn with probability   $\mu $, after an exponential holding time at $0$ each time it dies.)   We denote this rebirthed process by $\wt   X\!=\!
(\wt \Om, \wt  \FF_{t},  \wt  X_t, \wt \th_{t},\wt P^x)$.   To emphasize  the parameters that determine $\wt X$ we sometimes write $\wt X= \wt X(v^{\la}, m ,\mu,\bb)$.

 \bt\label{theo-mborel} 
The process  
 $  \wt X=  \wt X(v^{\la}, m ,\mu,\bb)$ is a recurrent Borel right process with state space $ S\cup \{0\}$ and $\la$--potential operator
 \begin{equation}
U^{\la}h(x)=V^{\la}h(x) +E^{x}\( e^{-\la \ze }\) \frac{h(0)+\bb \mu V^{\la}h}{\la+\bb  E^{\mu}\(1-e^{-\la \ze}   \)},   \label{193.3res}
 \end{equation}
and $\la$--potential densities, 
\begin{equation}
u^{\la}(x,y)= v^{\la}(x,y)+E^{x}\( e^{-\la \ze }\)\frac{1_{\{y=0\}}+\bb \mu v^{\la}(y)}{\la+\bb  E^{\mu}\(1-e^{-\la \ze} \)}\label{193.3a},
\end{equation}
with respect to   $\de_{0}(\cdot)+m$, where $\de_{0}(\cdot)$ is the unit mass at $0$
\et

 In Section \ref{sec-FBM} we take $S=(0,\ff)$, $\ze=\si$, $\mu=\bar \eta_{\ep}/\|  \bar \eta_{\ep}   \|$, $\bb=\| \bar  \eta_{\ep}   \|/p_{\ep}$ where $\lim_{\ep\to 0}p_{\ep}$ exists in $[0,\ff)$, so that by (\ref{6.2mm1}) we have
  \begin{equation}
U_{\ep}^{\la}h(x)=:  V^{\la}h(x) +E^{x}\( e^{-\la \si }\)\frac{p_{\ep} h(0) +  \bar \eta_{\ep} V^{\la}h} {   \la p_{\ep}+  E^{\bar \eta_{\ep}}\(1-e^{-\la \si} \)  },   \label{blu.1}
 \end{equation}
See the beginning of Section \ref{sec-FBM} for details. 
\medskip

  \noindent \textbf{Proof of Theorem \ref{theo-mborel}}       It is a classical result that   $ \wt X$  is a recurrent Borel right process. (See, e.g., \cite[Theorem 1]{Meyer} or \cite[(14.17)]{S}.) We continue to show that it has the $\la$--potential operator  (\ref{193.3res}).

Let   $t_n=\ze_0+\sum_{j=1}^{n}(J_{j}+\ze_j)$, where $\{ J_{j}, \ze_j \}_{j=1}^{\ff}$ are the consecutive     holding times and lifetimes of the rebirthed process, so that $ \ze_{n+1}= \ze\circ \wt\th_{t_n+J_{n+1}}$.

   We now take $h$ be bounded and consider  $U^\la $,   the $\la $--th potential operator for $\wt X$.
 For any $x\in S$,  \bea \label{6.6mk}
 && U^{\la }h(x)= \wt E^{x}\( \int_{0}^{\ff} e^{-\la t}h\(\wt X_{t}\) \,dt \) \\
  &&\quad=\nn  \wt E^{x}\(\int_{0}^{\ze_0}e^{-\la t}h\(\wt X_{t}\)\,dt\)+ \sum_{n=0}^{\ff}\wt E^{x}\(\int_{t_n}^{t_{n+1}}e^{-\la t}h\(\wt X_{t}\)\,dt\)\nn
\eea
and
\begin{eqnarray}
&&\wt E^{x}\(\int_{t_n}^{t_{n+1}}e^{-\la t}h\(\wt X_{t}\)\,dt\)
\label{7.12pre}
\\
&&=h(0)\,\,\wt E^{x}\(\int_{t_n}^{t_{n}+J_{n+1}}e^{-\la t}\,dt\)+\wt E^{x}\(\int_{t_n+J_{n+1}}^{t_{n+1}}e^{-\la t}h\(\wt X_{t}\)\,dt\)
\nonumber\\
&&=h(0)\,\,  \wt E^{x}\(    \frac{e^{-\la t_{n}}(1-e^{-\la J_{n+1}})}{\la}       \)\nonumber\\
&&+\wt E^{x}\(e^{-\la (t_n+J_{n+1})} \int_{0}^{ \ze_{n+1}}e^{-\la t}h\(\wt X_{t}\circ \wt\th_{t_n+J_{n+1}}\)\,dt\)\nn
\end{eqnarray}
 Using the fact that  $\wt X$ is a Markov process we have,
\bea \label{7.12}
&&\wt E^{x}\(e^{-\la (t_n+J_{n+1})} \int_{0}^{ \ze_{n+1}}e^{-\la t}h\(\wt X_{t}\circ\wt\th_{t_n+J_{n+1}}\)\,dt\)\\
&&\qquad =\wt E^{x}\(e^{-\la (t_n+J_{n+1})}\( \int_{0}^{ \ze }e^{-\la t}h\(X_{t}\)\,dt\)\circ\wt\th_{t_n+J_{n+1}}  \)\nn\\
&&
\qquad\wt E^{x}\(e^{-\la (t_n+J_{n+1})}\wt E^{\wt X_{t_n+J_{n+1}}}\( \int_{0}^{ \ze }e^{-\la t}h\(X_{t}\)\,dt\) \)\nn\\
&&\qquad
=\wt E^{x}\(e^{-\la (t_n+J_{n+1})} \) E^{\mu}\( \int_{0}^{ \ze }e^{-\la t}h\(X_{t}\)\,dt\),\nn
\eea
where we use  the fact that $\wt X_{t_n+J_{n+1}}$ has   probability distribution $\mu  $.    
Hence from (\ref{7.12pre})
\begin{eqnarray}
&&\wt E^{x}\(\int_{t_n}^{t_{n+1}}e^{-\la t}h\(\wt X_{t}\)\,dt\)
\label{7.12ps}
\\
&&=h(0)\,\,  \wt E^{x}\(    \frac{e^{-\la t_{n}}(1-e^{-\la J_{n+1}})}{\la}       \)+\wt E^{x}\(e^{-\la (t_n+J_{n+1})} \) E^{\mu}\( \int_{0}^{ \ze }e^{-\la t}h\(X_{t}\)\,dt\).\nn
\end{eqnarray}
Using the Markov property and  induction, in a similar manner  we see  that for $n\ge 0$,  \bea  \label{7.14qq}
&& \wt   E^x(e^{-\la (t_n+J_{n+1})}) = E^{x}\( e^{-\la \ze  }\)\wt   E^{0}\( e^{-\la (\sum^{n}_{i=1} \ze_{i}+\sum^{n+1}_{i=1}J_{i})        }\)\nn\\
&&
  =\(\frac{\bb}{\la +\bb}\)^{n+1} E^{x}\( e^{-\la \ze }\) \(  E^{\mu}\( e^{-\la  \ze }\)  \)^{n},  
\eea 
and
\begin{equation}
\wt   E^x(e^{-\la  t_n }) =\(\frac{\bb}{\la +\bb}\)^{n} E^{x}\( e^{-\la \ze }\) \(  E^{\mu}\( e^{-\la  \ze }\)  \)^{n}.\label{}
\end{equation}
It follows that
\begin{eqnarray}
&&\wt E^{x}\(    \frac{e^{-\la t_{n}}(1-e^{-\la J_{n+1}})}{\la }       \)\label{}
\\
&&=\(    \frac{(1- \(\frac{\bb}{\la +\bb}\))}{\la }       \)\(\frac{\bb}{\la +\bb}\)^{n} E^{x}\( e^{-\la \ze }\) \(  E^{\mu}\( e^{-\la  \ze }\)  \)^{n}
\nn
\\
&&=\(    \frac{1}{\la +\bb}       \)\(\frac{\bb}{\la +\bb}\)^{n} E^{x}\( e^{-\la \ze }\) \(  E^{\mu}\( e^{-\la  \ze }\)  \)^{n}.
\nonumber
\end{eqnarray} 
Thus by (\ref{7.12ps})
\begin{eqnarray}
&&\wt E^{x}\(\int_{t_n}^{t_{n+1}}e^{-\la t}h\(\wt X_{t}\)\,dt\)
\label{7.12pt}
\\
&&=\( h(0)\,\, \(    \frac{1}{\la +\bb}       \)+  \(\frac{\bb}{\la +\bb}\)\mu V^{\la }h\)\nn\\
&&\hspace{1 in}\times \(\frac{\bb}{\la +\bb}\)^{n} E^{x}\( e^{-\la \ze }\) \(  E^{\mu}\( e^{-\la  \ze }\)  \)^{n}.\nn
\end{eqnarray}
Note that  by simple integration and (\ref{6.2mm}),  
  \begin{equation}
E^{x}\( e^{-\la \ze}\)=1-\la E^{x}\( \int_{0}^{\ze}e^{-\la t} \,dt \)=1-\la V^{\la}1(x),  \label{rp.4}
\end{equation}
and   
\begin{equation} 
  E^{\mu}\( e^{-\la \ze}\)= 1-\la \mu V^{\la}1,  \label{rp.4a}
\end{equation}
which implies, in particular, that $\la \mu V^{p}1<1$.  
Thus by (\ref{7.12pt})
\begin{eqnarray}
&&\wt E^{x}\(\int_{t_n}^{t_{n+1}}e^{-\la t}h\(\wt X_{t}\)\,dt\)
\label{7.12pu}
\\
&&=E^{x}\( e^{-\la \ze }\)\( h(0)\,\, \(    \frac{1}{\la +\bb}       \)+  \(\frac{\bb}{\la +\bb}\)\mu V^{\la }h\)\nn\\
&&\hspace{1 in}\times      \(\(\frac{\bb}{\la +\bb}\)  \( 1-\la \mu V^{\la }1 \)\)^{n},\nn
\end{eqnarray}
so that
\begin{eqnarray}
&&\sum_{n=0}^{\ff}\wt E^{x}\(\int_{t_n}^{t_{n+1}}e^{-\la t}h\(\wt X_{t}\)\,dt\)
\label{7.12pv}
\\
&&=E^{x}\( e^{-\la \ze }\)\( h(0)\,\, \(    \frac{1}{\la +\bb}       \)+  \(\frac{\bb}{\la +\bb}\)\mu V^{\la }h\)\nn\\
&&\hspace{1 in}\times    \frac{1}{1-\(\(\frac{\bb}{\la +\bb}\)  \( 1-\la \mu V^{\la }1 \)\)}\nn\\
&&=E^{x}\( e^{-\la \ze }\)   \frac{\( h(0) +  \bb \mu V^{\la }h\) }{\la   \( 1+\bb\mu V^{\la }1 \)}. \nn
\end{eqnarray}
Consequently, using (\ref{6.6mk}) and (\ref{6.5mmma})
\begin{equation}
U^{\la }h(x)=V^{\la}h(x)+E^{x}\( e^{-\la \ze }\) \frac{\( h(0) +  \bb \mu V^{\la }h\)}{ \la+\bb  E^{\mu}\(1-e^{-\la \ze}   \)} . \label{rp.5}
\end{equation}
Since this is true for all $h$  we see that (\ref{193.3res}) follows. 
 \qed

 We assume that $X$ has continuous local times $ L_{t}^{y}$ normalized so that 
  \begin{equation}
  E^{x}\( L_{\ff}^{y}\)= v^{0}(x,y),\quad\,x,y\in S.
  \end{equation} 
  As in \cite[Lemma 2.1]{FMR} we can show that $\wt X$ has a jointly continuous local time $\{\wt L_{t}^{y}, y\in S, t\in R_{+}\}$  such that for each $\la>0$,
\bea
 && \wt E^{ x}\(\int_{0}^{\ff}  e^{-\la t}\, d_{t}\wt L^{y}_t\)=   u^{\la}(x,y)\label{int.3}\\
 &&\hspace{1 in}=v^{\la}(x,y)+E^{x}\( e^{-\la \ze }\)\frac{\bb \mu v^{\la}(y)}{\la+\bb  E^{\mu}\(1-e^{-\la \ze} \)} ,\quad\,x,y\in S. \nn
\eea
  
  We can prove an exact uniform modulus of continuity for the local times of 
$\wt X\!=\!
(\wt \Om, \wt  \FF_{t}, \wt  X_t, \wt \th_{t},\wt P^x)$ just as we did in \cite{FMR}. The only differences is that now we have holding times  in state $0$. As we shall see, the exact uniform modulus of continuity for the local times depends only on $X$, and not on the holding   rates nor the measure $\mu$.

 As before, we set $t_n=\ze_0+\sum_{j=1}^{n}(J_{j}+\ze_j)$, where $\{ J_{j}, \ze_j \}_{j=1}^{\ff}$ are the consecutive     holding times and lifetimes of the rebirthed process, so that $ \ze_{n+1}= \ze\circ \wt\th_{t_n+J_{n+1}}$. Set $t_{n}^{-}=\ze_0+\sum_{j=1}^{n-1}(J_{j}+\ze_j)+J_{n}$. We have $\wt \Om=(\Om \times R_{+}^{1})^{ \mathbf N}$ with elements $\wt\om=( \om_{0}, p_1, \om_{1}, p_2, \om_{2}, \ldots)$ and 
 $\wt P^{y}=P^{y}\times \prod_{n=1}^{\ff} \{P_{n, \bb} \times P_{n}^{\mu}\}$, where $P_{n, \bb}(dp_{i})$ are independent copies of the measure  $\bb e^{-\bb p}\,dp$ and $P_{n}^{\mu}$ are independent copies of $P^{\mu}$.

It follows as in the proof of \cite[Lemma 1.2]{FMR} that we have, 
  \begin{equation}
 \wt L^{x}_{t}(\wt\om)=L^{x}_{\ze(\om_{0})}(\om_{0})+\sum^{r-1}_{i=1}  L^{x}_{\ze(\om_{i})}(\om_{i}), \qquad x\in (0,1], \qquad t_{r-1}\leq t< t_{r}^{-}\label{k70.1}
 \end{equation} 
 and
   \begin{equation}
 \wt L^{x}_{t}(\wt\om)=L^{x}_{\ze(\om_{0})}(\om_{0})+\sum^{r-1}_{i=1}  L^{x}_{\ze(\om_{i})}(\om_{i})+   L^{x}_{t-t_{r}^{-}}(\om_{r}), \qquad x\in (0,1],\qquad t_{r}^{-}\leq t< t_{r}.\label{k70.1z}
 \end{equation}

  We now show that the terms on the right in (\ref{k70.1}) and (\ref{k70.1z}) are conditionally independent. Let $\ga,  \ga_{1}, \ga_{2},\ldots  $ denote  independent exponential random variables of mean $1/p$
 independent of everything else. For a probability measure $P$   we write 
 $P_{\ga}=P\times \,d\ga$ and similarly, 
 $P_{\ga_{i}}=P\times \,d\ga_{i}$.
 Let  $\ga' =\ga-t_{r}^{-}$. By (\ref{k70.1}) and (\ref{k70.1z})
 we have that
   \begin{equation}
 \wt L^{x}_{\ga}(\wt\om)=L^{x}_{\ze(\om_{0})}(\om_{0})+\sum^{r-1}_{i=1}  L^{x}_{\ze(\om_{i})}(\om_{i}), \qquad x\in (0,1], \qquad t_{r-1}\leq \ga< t_{r}^{-},\label{m70.2}
 \end{equation} 
 and
   \begin{equation}
 \wt L^{x}_{\ga}(\wt\om)=L^{x}_{\ze(\om_{0})}(\om_{0})+\sum^{r-1}_{i=1}  L^{x}_{\ze(\om_{i})}(\om_{i})+   L^{x}_{\ga'}(\om_{r}), \qquad x\in (0,1],\qquad t_{r}^{-}\leq \ga< t_{r}.\label{m70.2z}
 \end{equation}

  Let
 \bea
  \mathcal{L}_{0}&=&\{  L^{x}_{\ze(\om_{0})}(\om_{0}) , x\in (0,1]\},   \label{}\\
 \mathcal{L}_{i}&=&\{    L^{x}_{\ze(\om_{i})}(\om_{i}), x\in (0,1]\},\qquad 1\le i\leq r-1, \nn  \\
 \mathcal{L}'_{r}&=&  \{  L^{x}_{\ga'\wedge \ze(\om_{r})}(\om_{r}), x\in (0,1]\},\nn\\
   \ov {\mathcal{L}} _{r} &=&\nn \{L^{x}_{\ga_{r} }(\om_{r}), x\in (0,1]\}.\label{}\eea
In this notation (\ref{m70.2}) and (\ref{m70.2z}) can be written as,
     \begin{equation}
 \wt L^{x}_{\ga}(\wt\om)=\sum^{r-1}_{i=0}   \mathcal{L}_{i},  \qquad t_{r-1}\leq \ga< t_{r}^{-}\label{m70.2a}
 \end{equation}
 and
      \begin{equation}
 \wt L^{x}_{\ga}(\wt\om)=\sum^{r-1}_{i=0}   \mathcal{L}_{i} +  \mathcal{L}'_{r},\qquad  t_{r}^{-}\leq \ga<t_{r}.\label{m70.2az}
 \end{equation}

\bl\label{lem-mcondind}
 For all $y>0$, $r\geq 2$ and all Borel sets, $B_0, B_1, B_2,\ldots, B_r$, in $C((0,1], R^1)$,
  \begin{eqnarray}
 &&\wt P_{\ga}^{y}\( \mathcal{L}_{0}\in B_0, \mathcal{L}_{1}\in B_1, \ldots,  \mathcal{L}_{r-1}\in B_{r-1}\,\Big |\, t_{r-1}\leq \ga< t_{r}^{-} \)
\nn
 \\
 &&= P_{\ga_{0}}^{y}\( \mathcal{L}_{0}\in B_0 \,\Big |\, \ze(\om_{0})<\ga_{0}   \) \prod^{r-1}_{i=1}   P^{\mu}_{\ga_{i}} \(\mathcal{L}_{i}\in B_i\,\Big |\,  \ze(\om_{i})<\ga_{i} \),
 \label{m70.5mp}
 \end{eqnarray}  
 and
   \begin{eqnarray}
 &&\wt P_{\ga}^{y}\( \mathcal{L}_{0}\in B_0, \mathcal{L}_{1}\in B_1, \ldots,  \mathcal{L}_{r-1}\in B_{r-1},  \mathcal{L}'_{r} \in B_r \,\Big |\,  t_{r}^{-}\leq \ga< t_{r} \)
\nn
 \\
 &&\qquad= P_{\ga_{0}}^{y}\( \mathcal{L}_{0}\in B_0 \,\Big |\, \ze(\om_{0})<\ga_{0}   \)  \prod^{r-1}_{i=1}  P^{\mu}_{\ga_{i}} \(\mathcal{L}_{i}\in B_i\,\Big |\, \ze(\om_{i})<\ga_{i}\)\nn\\
 &&\hspace{1 in}  \times    P^{\mu}_{\ga_{r}} \( \ov {\mathcal{L}}_{r}\in B_r\,\Big |\,  \ga_{r}\leq \ze(\om_{r})\).
 \label{m70.5mpa}
 \end{eqnarray} 
 In particular,
 \begin{equation} \label{m3.26}
   \wt P_{\ga}^{y}\(  \mathcal{L} '_{r}\in B_r\,\Big |\,  t_{r}^{-}\leq \ga< t_{r}\)=    P^{\mu}_{\ga_{r}} \( \ov {\mathcal{L}}_{r}\in B_r\,\Big |\,  \ga_{r}\leq \ze(\om_{r})\),
\end{equation}
and (\ref{m70.5mpa}) implies that for each $r\geq 2$,  
 \begin{equation}
 \mathcal{L}_{0}, \, \mathcal{L}_{1}, \,\ldots, \, \mathcal{L}_{r-1},\,  \mathcal{L}'_{r}, \label{m70.0}
 \end{equation}
 are conditionally independent, given that $ t_{r}^{-}\leq \ga< t_{r}$. 

  \el
  
 \medskip
   \Proof Let $\mathcal{C}$ be a dense countable subset of $(0,1]$.    For any function $f(x)$,   $x\in\CC$,    set  
\begin{equation} \label{}
   (b_{i},f)=\sum_{x\in \mathcal{C}}a_{i, x} f(x),
\end{equation}  where the $a_{i, x}\in R^{1}$, $0\le i\le r$    and $\sum_{x\in \mathcal{C}}|a_{i, x}|<\ff$  for each $i$. 
We begin this proof by calculating the Laplace transform of,  
\begin{equation} \label{}
(b_{0}, L^{\cdot}_{\ze(\om_{0})}(\om_{0})) +  \sum^{r-1}_{i=1} (b_{i},   L^{\cdot}_{\ze(\om_{i})}(\om_{i}))+  (b_{r},   L^{\cdot}_{\ga'}(\om_{r})).
\end{equation} 
We have,
 \begin{eqnarray}
 \lefteqn{\wt E_{\ga}^{y}\(\exp\({- (b_{0}, L^{\cdot}_{\ze(\om_{0})}(\om_{0})) -   \sum^{r-1}_{i=1} (b_{i},   L^{\cdot}_{\ze(\om_{i})}(\om_{i}))-  (b_{r},   L^{\cdot}_{\ga'}(\om_{r}))}\);  t_{r}^{-}\leq \ga< t_{r} \)}
\nn
 \\
 &&=\wt E^{y} \int_{ t_{r}^{-}}^{t_{r}} pe^{-pt}\exp\({-(b_{0}, L^{\cdot}_{\ze(\om_{0})}(\om_{0})) - \sum^{r-1}_{i=1} (b_{i},   L^{\cdot}_{\ze(\om_{i})}(\om_{i}))-  (b_{r},   L^{\cdot}_{t-t_{r}^{-} }(\om_{r})}\) \,dt 
 \nonumber \\
 &&=\wt E^{y}\int_{ 0}^{ t_{r}-t_{r}^{-}} pe^{-p(u+ t_{r}^{-})} \label{m70.3}\\
 && \hspace{1 in} \exp\({- (b_{0}, L^{\cdot}_{\ze(\om_{0})}(\om_{0})) -\sum^{r-1}_{i=1} (b_{i},   L^{\cdot}_{\ze(\om_{i})}(\om_{i}))-  (b_{r},   L^{\cdot}_{u}  (\om_{r}))}\)  \,du.\nn    
 \eea
 Since, obviously, for any function $f(u)$
\begin{equation} \label{}
\int_{ 0}^{ t_{r}-t_{r}^{-}} pe^{-p(u+ t_{r}^{-})} f(u) \,du=e^{-p t_{r}^{-} } \int_{ 0}^{ t_{r}-t_{r}^{-}} pe^{-p u }f(u) \,du,
\end{equation}
we can write the last two lines of (\ref{m70.3}) in the form,
\be 
I= \wt E^{y}\( e^{-p t_{r}^{-} }e^{-(b_{0}, L^{\cdot}_{\ze(\om_{0})}(\om_{0}))- \sum^{r-1}_{i=1} (b_{i},   L^{\cdot}_{\ze(\om_{i})}(\om_{i}))}    \int_{ 0}^{\ze_{r}} pe^{-pu} e^{-    (b_{r},   L^{\cdot}_{u}  (\om_{r}))}   \,du\).
  \ee
Since $t_{r}^{-}=\ze_0+\sum_{j=1}^{r-1}(J_{j}+\ze_j)+J_{r}$
 and  the holding times $J_{i}$ are independent, 
 \begin{equation}
\wt E^{y}\(  e^{-p\(\sum_{i=1}^{r}J_{i}\)}  \)=\prod_{i=1}^{r}\wt E^{y}\(  e^{-p J_{i} }  \) = \(\frac{\bb}{\bb+p}\)^{r}. \label{}
 \end{equation}
In addition, since the holding times $J_{i}$ are independent of everything else,
  \begin{eqnarray}
  &&I=\(\frac{\bb}{\bb+p}\)^{r} \wt E^{y}\( e^{-p \ze(\om_{0})}e^{-(b_{0}, L^{\cdot}_{\ze(\om_{0})}(\om_{0}))} \prod_{i=1}^{r-1} e^{-p \ze(\om_{i})}e^{-(b_{i}, L^{\cdot}_{\ze(\om_{i})}(\om_{i}))}\right.\nn\\
  &&\left. \hspace{2.2 in} \int_{ 0}^{\ze(\om_{r})} pe^{-pu} e^{-    (b_{r},   L^{\cdot}_{u}  (\om_{r}))}   \,du\).
  \label{}
  \end{eqnarray}
Using the Markov property, we then see that
  \begin{eqnarray}
  &&I=\(\frac{\bb}{\bb+p}\)^{r}E^{y}\(e^{-p \ze(\om_{0})} e^{-(b_{0}, L^{\cdot}_{\ze(\om_{0})}(\om_{0}))}\) \prod_{i=1}^{r-1} E^{\mu}\(e^{-p \ze(\om_{i})}e^{-(b_{i}, L^{\cdot}_{\ze(\om_{i})}(\om_{i}))}\)\nn\\
  && \hspace{2.6 in}\times E^{\mu}\(\int_{ 0}^{\ze(\om_{r})} pe^{-pu} e^{-    (b_{r},   L^{\cdot}_{u}  (\om_{r}))}   \,du\)\nn\\
 &&= \(\frac{\bb}{\bb+p}\)^{r} \,\,E_{\ga_{0}}^{y}\(e^{-(b_{0}, L^{\cdot}_{\ze(\om_{0})}(\om_{0}))} ; \ze(\om_{0})<\ga_{0} \)  \prod^{r-1}_{i=1}    E^{\mu}_{\ga_{i}}\(e^{-(b_{i}, L^{\cdot}_{\ze(\om_{i})}(\om_{i}))} ; \ze(\om_{i})<\ga_{i} \) \nonumber\\
 &&\hspace{2 in} \times   E^{\mu}_{\ga_{r}}\(e^{-  (b_{r},   L^{\cdot}_{\ga_{r}}(\om_{r}))} ; \ga_{r}\leq \ze(\om_{r})  \).\label{m70.3a}
 \end{eqnarray}
 Setting all the $b_{i}=0$ we see that,
 \be \wt P_{\la}^{y} \( t_{r}^{-}\leq \ga< t_{r} \)=\(\frac{\bb}{\bb+p}\)^{r} \,\,P_{\ga_{0}}^{y}\( \ze(\om_{0})<\ga_{0} \)\prod^{r-1}_{i=1}  E^{\mu}_{\ga_{i}} \(\ze(\om_{i}<\ga_{i}\) E^{\mu}_{\ga_{r}} \(  \ga_{r}\leq \ze(\om_{r}\).
 \ee
 Consequently,
 \begin{eqnarray}
&&\wt E_{\ga}^{y}\(\exp\({- (b_{0}, L^{\cdot}_{\ze(\om_{0})}(\om_{0})) -   \sum^{r-1}_{i=1} (b_{i},   L^{\cdot}_{\ze(\om_{i})}(\om_{i}))-  (b_{r},   L^{\cdot}_{\ga'}(\om_{r}))}\)\,|\,  t_{r}^{-}\leq \ga<\leq t_{r} \)
\nn
 \\
 &&=  E_{\ga_{0}}^{y}\(e^{-(b_{0}, L^{\cdot}_{\ze(\om_{0})}(\om_{0}))} \,|\, \ze(\om_{0})<\ga_{0} \)  \prod^{r-1}_{i=1}   E^{\mu}_{\ga_{i}}\(e^{-(b_{i}, L^{\cdot}_{\ze(\om_{i})}(\om_{i}))}  \,|\,  \ze(\om_{i})<\ga_{i} \)  \label{70.5}\\
 &&\hspace{1 in} \times       E^{\mu}_{\ga_{r}}\(e^{-  (b_{r},   L^{\cdot}_{\ga_{r}}(\om_{r}))} \,|\,\ga_{r}\leq \ze(\om_{r})  \).
 \nonumber
 \end{eqnarray}
 
 Taking the inverse  Laplace transform    gives  (\ref{m70.5mpa}) with the $B_{i}$  restricted to $\mathcal{C}$. The statement in (\ref{m70.0}) follows because the local times are continuous.    (\ref{m70.5mp})  follows similarly.    \qed

 Let 
 $ L^{x}_{i, t}$    be independent copies of $L^{x}_{t}$  and  let $\ze(i)$ denote the death time of  $L^{x}_{i,t}$, $i=1, 2, \ldots  $. 
 For any set $C$, let $F(C)$ denote the set of real--valued functions
$f$ on $C$. Define the evaluations
$\it{i}_{ x}:F(C)\mapsto R^{{\rm 1}}$ by $\it{i}_{ x}( f)=f( x)$. We use
$\mathcal{M}( F(C))$ to denote  the smallest $\si$-algebra for which the
evaluations
$\it{i}_{ x}$ are Borel measurable for all
$x\in C$. ($\mathcal{M}( F(C))$ is generally referred to  as the
$\si$-algebra  of cylinder sets in $F(C)$.)

 \bt\label{theo-m1.1} When, for some $r\geq 1$, and measurable set of functions $B\in \mathcal{M}( F((0,1]))$,
 \begin{equation}
\(P_{0}^{y}\times \prod^{r-1}_{i=1} P^{\mu}_{i}\times   P^{\mu}_{r, \ga_{r}}\)  \( \sum^{r-1}_{i=0} L^{\cdot}_{i, \ff}+  L^{\cdot}_{r,\ga_{r}}\in B\)=1,  \label{m70.6}
 \end{equation}
 then  
  \begin{equation}
\wt P_{\ga}^{y} \(\wt L^{\cdot}_{\ga}(\wt\om)\in B\,\Big |\, t_{r}^{-}\leq \ga< t_{r}\)=1,\label{m70.7}
 \end{equation}
 with a similar result for $t_{r-1}\leq \ga< t_{r}^{-}$.
 \et
 
 \Proof  
 Since $\{\ga_{i}, 0\le i\le r-1\}$ do not appear in the event  $\{ \sum^{r-1}_{i=0} L^{\cdot}_{i, \ff}+  L^{\cdot}_{r,\ga_{r}}\in B\}$, we can  write  (\ref{m70.6}) as,
  \begin{equation}
\(P_{0,\ga_{0}}^{y}\times \prod^{r-1}_{i=1} P^{\mu}_{i,\ga_{i}} \times   P^{\mu}_{r, \ga_{r}} \)  \( \sum^{r-1}_{i=0} L^{\cdot}_{i, \ff}+  L^{\cdot}_{r,\ga_{r}}\in B\)=1. \label{m70.6c}
 \end{equation}
Now using the fact that  an almost sure event   with respect to a given probability occurs almost surely for any conditional version of the probability, we see 
  that   (\ref{m70.6c}) implies that,     
  \be
\(P_{0,\ga_{0}}^{y}\times \prod^{r-1}_{i=1} P^{\mu}_{i,\ga_{i}} \times   P^{\mu}_{r, \ga_{r}} \)\( \sum^{r-1}_{i=0} L^{\cdot}_{i, \ff}+  L^{\cdot}_{r,\ga_{r}}\in B\,\Big |\,\mathcal{A}\)=1,
  \label{m70.6d}
\ee 
 where
  \begin{equation}
 \mathcal{A}=\{ \ze(0)<\ga_{0}, \, \ze(1)<\ga_{1},       \,\cdots,\,  \ze(r-1)<\ga_{r-1},  \,  \ga_{r}\leq \ze(r)\}.    \label{}
  \end{equation}
 We note that   for each $0\leq i\leq r-1$, $L^{x}_{\ze(\om_{i})}(\om_{i})$ has the law of $L^{x}_{\ff}$,
 and $ \ov {\mathcal{L}}= L^{x}_{\ga_{r}}(\om_{r})$  has the law of $ L^{x}_{\ga}$. Hence (\ref{m70.6d}) is equivalent to 
    \begin{eqnarray}
  &&\(P_{0,\ga_{0}}^{y}\times \prod^{r-1}_{i=1} P^{\mu}_{i,\ga_{i}} \times   P^{\mu}_{r, \ga_{r}} \)
  \label{m70.6e}
  \\
  &&\qquad  \(   \sum^{r-1}_{i=0} \mathcal{L}_{i}+\ov {\mathcal{L}}\in B\,\Big |\,    \ze(\om_{0})<\ga_{0},\cdots,  \ze(\om_{r-1})<\ga_{r-1},  \,  \ga_{r}\leq \ze (\om_{r}) \)=1.
  \nonumber
  \end{eqnarray}
It then follows from (\ref{m70.5mpa}) that, 
  \begin{equation}
  \wt P_{\ga}^{y}\(\sum^{r-1}_{i=0} \mathcal{L}_{i}+\mathcal{L}'_{r}\in B   \,\Big |\,     t_{r}^{-}\leq \ga< t_{r} \)=1.\label{}
  \end{equation}
Using (\ref{m70.2az})  we see that this is the statement in (\ref{m70.7}).
\qed

Set $S=(0,1]$.
\bt\label{theo-ITcond}  Let $\{\eta_{i, 0} (x)$, i=1,\ldots\} be independent  Gaussian processes with covariance $ v^{0}(x,y)$. Let   $\eta_{p} (x)$ be a Gaussian process with covariance $  v^{p}(x,y)$ independent of the $\{\eta_{i, 0} (x)\}$. For any $r\geq 1$ set $G_{r, s}=\{G_{r, s}(x),x\in S \}$ where
\be G_{r, s}\(x\)=\sum_{i=0}^{r-1} \frac{1}{2}(\eta_{i, 0 }(x)+s)^2+\frac{1}{2}(\eta_{p}(x)+s)^2, \label{m2.21p}
\ee
and $ P_{G_{r, s}} =  \prod_{i=0}^{r-1}  P_{\eta_{i, 0} } \times   P_{\eta_{p}} $. When 
\begin{equation}
 P_{G_{r, s}} \(G_{r, s}\in B\)=1,  \label{m70.6n}
 \end{equation}
  for a measurable set of functions $B\in \mathcal{M}( F(S))$, then 
  \begin{equation}
\(\wt P_{\ga}^{y}\times P_{G_{r, s}}\) \(\wt L^{x}_{\ga}(\wt\om)+G_{r, s} \in B\,\Big |\,    t_{r}^{-}\leq \ga< t_{r}\)=1,\qquad   \forall y>0.\label{70.7n}
 \end{equation}
\et

\noindent\textbf{Proof of Theorem \ref{theo-ITcond} }
 It follows from the Eisenbaum Isomorphism Theorem, \cite[Theorem 8.1.1]{book}, that,   
\bea\lefteqn{
\Big\{ L^x_{i,\ff}+\frac{1}{2}(\eta_{i, 0} (x)+s)^2\,\,;\,x\in S\,,\,P_{i}^y\times
P_{\eta_{i, 0}}\Big\}\label{mit1.2f}}\\ &&\qquad\stackrel{law}{=}
\Big\{\frac{1}{2}(\eta_{i, 0} (x)+s)^2\,\,;\,x\in S\,,\,\(1+\frac{\eta_{0,i} (y)}s\)P_{\eta_{i, 0}}\Big\},\nn
\eea
where $P_{\eta_{i, 0}}$ is the probability  of $\eta_{i, 0}$ and $P_{i}^y$ the probability of $X$ started at $y$. 
It also follows from the Eisenbaum Isomorphism Theorem that,
\bea\lefteqn{
\Big\{ L^x_{r, \ga_{r}}+\frac{1}{2}(\eta_{p} (x)+s)^2\,\,;\,x\in S\,,\,P_{r, \ga_{r}}^y\times
P_{\eta_{p}}\Big\}\label{mit1.2}}\\ &&\qquad\stackrel{law}{=}
\Big\{\frac{1}{2}(\eta_{p} (x)+s)^2\,\,;\,x\in S\,,\,(1+\frac{\eta_{p } (y)}s)P_{\eta_{p}}\Big\}\nn,
\eea
where $P_{\eta_{p}}$ is the probability of $\eta_{p}$  and $P_{r, \ga_{r}}^y= P_{r}^y\times d\ga_{r}$. 

Set
\begin{equation} \label{m3.37} 
  \eta_{i, 0} (\mu)=\int_S \eta_{i, 0} (x)\,d\mu(x),\qquad\text{and}\qquad \eta_{p} (\mu)=\int_S \eta_{p} (x)\,d\mu(x).
\end{equation} 
It follows from (\ref{mit1.2f}) and (\ref{mit1.2})   that  for each  $r\geq 1$,
\bea
&& 
\Big\{ \sum^{r-1}_{i=0} L^{x}_{i, \ff}+ L^{x}_{r, \ga_{r}}+\sum_{i=0}^{r-1}\frac{1}{2}(\eta_{i, 0}(x)+s)^2+\frac{1}{2}(\eta_{p}(x)+s)^2\,;\,\, x\in S,\,\,\nn  \\
&&\hspace{1 in}  \( P_{0}^{y}\times \prod^{r-1}_{i=1} P^{\mu}_{i}\times   P^{\mu}_{r, \ga_{r}} \)     \times\prod_{i=0}^{r-1}  P_{\eta_{i, 0} } \times   P_{\eta_{p} } \Big\} \label{m3.38}\\
 &&\stackrel{law}{=}
\Big\{ \sum_{i=0}^{r-1}\frac{1}{2}(\eta_{i, 0}(x)+s)^2+\frac{1}{2}(\eta_{p}(x)+s)^2\,\,\,\,\,;\,x\in S\,,\nn\\
&& \hspace{.5 in}\qquad     \,(1+\frac{\eta_{ 0,0}(y)}s)P_{\eta_{0,0} } \times\prod_{i=1}^{r-1}  (1+\frac{\eta_{i, 0}(\mu)}s)P_{\eta_{i, 0} }\times   ( (1+\frac{\eta_{p}(\mu)}s)P_{\eta_{p} }) \Big\}.\nn
\eea
  For $r\geq 1$ it follows from (\ref{m70.6n}) and (\ref{m3.38}) that 
 \be
 \(P_{0}^{y}\times \prod^{r-1}_{i=1} P^{\mu}_{i } \times   P^{\mu}_{r, \ga_{r}} \times   P_{G_{r, s} }\)
\( \sum^{r-1}_{i=0} L^{\cdot}_{i, \ff}+ L^{\cdot}_{r,\ga_{r}}+G_{r,s}(\cdot)\in B\)=1. \label{}
\ee
Therefore, for $P_{G_{r, s} }$ almost every $\om'$,
 \be
 \(P_{0}^{y}\times \prod^{r-1}_{i=1} P^{\mu}_{i } \times   P^{\mu}_{r, \ga_{r}} \)
\( \sum^{r-1}_{i=0} L^{\cdot}_{i, \ff}+ L^{\cdot}_{r,\ga_{r}}+G_{r,s}(\cdot, \om')\in B\)=1. \label{}
\ee
 This  theorem now follows    from   Theorem \ref{theo-m1.1}.\qed

 \section{Exact uniform moduli of continuity for the local times of Feller Brownian motions}\label{sec-FBM}
 
   Let $S=(0,\ff)$ with the usual Borel $\si$-algebra.  
Let   $X=(\Om,  \FF_{t},  X_t,\th_{t},   P^x)$ be Brownian motion started at $x>0$ and killed at $\si$, the first time it reaches $0$. $Q_{t}$, the transition  semi-group of $X$, has density
\begin{equation}
q_{t}(x,y)=\frac{e^{ -(y-x)^{2}/2t}}{\sqrt{2\pi t}}-\frac{e^{ -(y+x)^{2}/2t}}{\sqrt{2\pi t}},\quad x,y>0\label{}
\end{equation}
 with respect to Lebesgue measure on $(0,\ff)$.
$X$ has continuous strictly positive  $\la-$potential   densities   with respect to Lebesgue measure on $(0,\ff)$:
\be
 v^{\la}(x,y)=  \frac{e^{-\sqrt{2\la}|y-x|}}{\sqrt{2\la}}-\frac{e^{-\sqrt{2\la}|y+x|}}{\sqrt{2\la}}, \quad x,y>0,\quad \la> 0,
\ee
and
\be
 v^{0}(x,y)=  2(|x|\wedge |y|), \quad x,y>0.
\ee
For measurable functions  $h$ on $(0,\ff)$ we  define  the $\la$--th potential operator,
\begin{equation}
V^{\la}h(x)=\int_{0}^{\ff}   v^{\la}(x,y)h(y) \,dy.    \label{su.1}
\end{equation}
 Using (\ref{6.2mm1}) we have 
 \begin{equation}
 \la\,V^{\la}1(x)=E^{ x }\(1-e^{-\la \si}   \)=1-e^{-x\sqrt{2\la}},\label{su.2}
 \end{equation}
for any $x>0$, see \cite[p. 56]{Blu} or \cite[(3.107)]{book}.

We note that reflecting Brownian motion in $[0,\ff)$ has continuous strictly positive  $\la-$potential   densities   with respect to Lebesgue measure on $[0,\ff)$:
\be
 r^{\la}(x,y)=  \frac{e^{-\sqrt{2\la}|y-x|}}{\sqrt{2\la}}+\frac{e^{-\sqrt{2\la}|y+x|}}{\sqrt{2\la}}, \quad x,y\geq 0. 
\ee 
For measurable functions  $h$ on $[0,\ff)$ we  define  the $\la$--th potential operator,
\begin{equation}
R^{\la}h(x)=\int_{0}^{\ff}   r^{\la}(x,y)h(y) \,dy.    \label{}
\end{equation}
It is easy to check that
 \begin{equation}
\lim_{x\to 0}x^{-1}V^{\la}h(x)=\sqrt{2\la}\,R^{\la}h(0). \label{blu.40r}
 \end{equation}

We now show how to obtain exact uniform moduli of continuity for the local times of Feller Brownian motions. A Feller Brownian motion $\wt X$ is a recurrent Borel right process with state space $[0,\ff)$ and $1$- potential operator 
  \begin{equation}
 U^{1}h(x)=V^{1}h(x) +e^{-x\sqrt{2 }}\,\,\,\(p h(0)+q R^{1}h(0)+  \eta V^{1}h\),   \label{blu.100}
 \end{equation}
for $p,q\geq 0$ and $\eta $ a measure on $S=(0,\ff)$ such that 
\begin{equation}
p+q+\int_{0}^{\ff}(1-e^{-x\sqrt{2}})\,\eta(dx)=1.\label{blu.101}
\end{equation}

When restricted to $(0,\ff)$,  $U^{1}$ has a continuous density $u ^1 (x,u)$ with respect to  Lebesgue measure. Let $\wt L^{u}_{t}, \,u\in S, \,t\in R_{+}$ denote the local times of $\wt X$, normalized so that for
 $y, u>0$,
 \begin{equation}
  \wt E^{ y}\(\int_{0}^{\ff}  e^{-  t}\, d_{t}\wt L^{u}_{t}\)=   u ^1 (y,u). \label{int.0}
\end{equation}
The local times $\wt L^{u}_{t} $ are jointly continuous in $u,t$.
Let $c_{j}\downarrow 0$, and let  $\mathcal{C}$ be a countable dense subset of $(0,1]$.
Using the joint continuity of $\wt L^{u}_{t} $, to prove Theorem \ref{theo-uniformmodulus} it suffices to show that for any   $j$ and $y>0$
   \be
\lim_{h\to 0}\sup_{\stackrel{|u-v|\le h }{ u,v\in [c_{j},1]\cap \mathcal{C} }} \frac{|\wt L_{t}^{u}-\wt L_{t}^{v}|}{  (|u-v| \log  1/| u-v|)^{1/2} }= \sup_{u\in [c_{j},1]\cap \mathcal{C}}\(2\wt L_{t}^{u}\)^{1/2},  \quad a.e.\,\,\, t,  \,\,\, \wt P^{y}\,\,
a.s.\label{int.1cut}
\ee


 \medskip Define the entrance law
 \begin{equation}
 \eta_{t}=\eta Q_{t},\label{elaw.1}
 \end{equation}
  so that 
 \begin{equation}
\eta_{s}Q_{t}= \eta_{s+t}, \quad \forall s,t>0.\label{elaw.2}
 \end{equation}
 By   (\ref{blu.101}) we have  $\eta V^{1}1=\int_{0}^{\ff}(1-e^{-x\sqrt{2}})\,\eta(dx)<\ff$
so that for bounded $h$
 \begin{equation}
\lim_{\ep\to 0} \,\,\eta_{\ep} V^{1}h=\eta V^{1}h<\ff.\label{blu.40}
 \end{equation}
Even when $\|\eta\|=\ff$,  the $\eta_{t}$'s are finite measures. See \cite[p. 138]{Blu}.

We will obtain exact uniform moduli of continuity for the local times of Feller Brownian motions by approximating them by the local times of the holding and jumping process described in (\ref{blu.1}) with $\la$-potential operator
  \begin{equation}
U_{\ep}^{\la}h(x)=: V^{\la}h(x) +E^{x}\( e^{-\la \si }\)\frac{p_{\ep} h(0) +  \bar \eta_{\ep} V^{\la}h} {   \la p_{\ep}+  E^{\bar \eta_{\ep}}\(1-e^{-\la \si} \)  }.   \label{blu.1p}
 \end{equation}
Here  $E^{x}\( e^{-\la \si }\) =e^{-\sqrt{2\la}\,\, x}$, and thus $ E^{\bar\eta_{\ep}}\(1-e^{-\la \si}   \)=\int_{0}^{\ff}\(1-e^{-\sqrt{2\la}\,\, x}\)\,d\bar\eta_{\ep}(x)$. Thus we can write 
(\ref{blu.1p})   as
  \begin{equation}
U_{\ep}^{\la}h(x)=V^{\la}h(x) +e^{-\sqrt{2\la}\,\, x}\,\,\frac{p_{\ep}h(0)+  \bar\eta_{\ep} V^{\la}h}{\la p_{\ep} +\int_{0}^{\ff}\(1-e^{-\sqrt{2\la}\,\, x}\)\,d\bar\eta_{\ep}(x)}.   \label{blu.17}
 \end{equation}
 In  (\ref{blu.17}) take  
 \begin{equation}
\bar\eta_{\ep}=\frac{q}{\sqrt{2}\,\ep\,}\,\de_{\ep}+\eta_{\ep}\label{blu.37}
\end{equation}
with $\eta_{\ep}$ as in (\ref{elaw.1}) to obtain  
  \begin{equation}
U_{\ep}^{\la}h(x)=V^{\la}h(x) +e^{-\sqrt{2\la}\,\, x}\,\,\frac{p_{\ep}h(0)+ \frac{q}{\sqrt{2}\,\ep\,}\, V^{\la}h(\ep) + \eta_{\ep} V^{\la}h}{\la p_{\ep}+q\frac{\(1-e^{-\sqrt{2\la }\,\,\ep}\)}{\sqrt{2}\,\ep\,} +\int_{0}^{\ff}\(1-e^{-\sqrt{2\la}\,\, x}\)\,d\eta_{\ep}(x)}.   \label{blu.17q}
 \end{equation}
With $\la=1$, 
 (\ref{blu.17q}) then becomes 
   \be
U_{\ep}^{1}h(x)=V^{1}h(x) +e^{-\sqrt{2 }\,\, x}\,\,\frac{p_{\ep}h(0)+ \frac{q}{\sqrt{2}\,\ep\,}\, V^{1}h(\ep) +\eta_{\ep} V^{1}h}{p_{\ep}+q\frac{\(1-e^{-\sqrt{2 }\,\,\ep}\)}{\sqrt{2}\,\ep\,}+\int_{0}^{\ff}\(1-e^{-\sqrt{2}\,\, x}\)\,d\eta_{\ep}(x)}.  \label{blu.41}
 \ee
 Restricted to  $(0,\ff)$, $U_{\ep}^{1}$ has continuous strictly positive  $1-$potential   densities  $u_{\ep} ^1 (y,u)$ with respect to Lebesgue measure.
 
 If the $p$ in (\ref{blu.100}) is strictly greater than zero, set the $p_{\ep}$ here equal to $p$. Then using (\ref{blu.40}), (\ref{blu.40r}) and (\ref{blu.101}) we see that
    \be
\lim_{\ep\to 0}U_{\ep}^{1}h(x)=V^{1}h(x) +e^{-\sqrt{2 }\,\, x}\,\,\( ph(0)+ q\, R^{1}h(0) +\eta V^{1}h\)= U^{1}h(x),  \label{blu.41f}
 \ee
and hence
\begin{equation}
\lim_{\ep\to 0} u_{\ep} ^1 (y,u)=u^1 (y,u),\quad y,u>0.\label{int.0dense}
\end{equation}

  If the $p$ in (\ref{blu.100}) is equal to zero, take   $p_{\ep}\to 0$. As before we find that $\lim_{\ep\to 0}U_{\ep}^{1}h(x)=U^{1}h(x) $ and $\lim_{\ep\to 0}u_{\ep} ^1 (y,u)=u^1 (y,u) $.

 We first consider the  case when $\lim_{\ep\to 0}\|\bar\eta_{\ep}\|=\ff.$ Let $ \mathcal{C}$ be a countable dense subset of $(0,1]$ such that 
 $\eta(\mathcal{C})=0$.   It follows from  (\ref{blu.101}) that $\eta$ is a countable sum of finite measures, so that we can find such a $\mathcal{C}$.  Since $\wt X$ has a countable number of excursions, it follows that almost surely none of them begin in $\mathcal{C}$.
 That is, almost surely 
 \begin{equation}
  Y_{s}(0)\notin \mathcal{C} \text{ for all } s \text{ such that }\tau(s^{-})<\tau(s),\label{lil.c0}
 \end{equation}
  using the notation of \cite[p. 145]{Blu}.
  We will prove (\ref{int.1cut}) by constructing a family of local times $\{ \wt L_{0,t}^{u}, \,\,\, \,u\in \mathcal{C}\}$
 such  that  
 \begin{equation}
  \wt E^{ y}\(\int_{0}^{\ff}  e^{-  t}\, d_{t}\wt L^{u}_{0,t}\)=   u^1 (y,u),\quad y>0, \label{int.0ep0}
\end{equation}
which implies by \cite[Theorem 3.6.7]{book} that 
\begin{equation}
  \wt L_{0,t}^{u}=\ \wt L_{t}^{u}, \,\,\, \forall t\geq 0, \,\,a.s. \label{lil.q1r}
\end{equation}
 (\ref{int.1cut}) will then follow once we will show that
 for any  sequence  $c_{j}\to \ff$ and $y>0$
   \be
\lim_{h\to 0}\sup_{\stackrel{|u-v|\le h }{ u,v\in [c_{j},1]\cap \mathcal{C} }} \frac{|\wt L_{0,t}^{u}-\wt L_{0,t}^{v}|}{  (|u-v| \log  1/| u-v|)^{1/2} }= \sup_{u\in [c_{j},1]\cap \mathcal{C}}\(2\wt L_{0,t}^{u}\)^{1/2},  \quad a.e.\,\,\, t,  \,\,\, \wt P^{y}\,\,
a.s.\label{int.1cutq}
\ee

Let $\wt L_{\ep,t}^{u}, \,u\in S, \,t\in R_{+}$ denote the local times of the process  
 $\wt   X_{\ep,t}=\wt X(v^{p}, m ,\mu_{\ep},\bb_{\ep})$, normalized so that for
 $y, u>0$,
 \begin{equation}
  \wt E^{ y}\(\int_{0}^{\ff}  e^{-  t}\, d_{t}\wt L^{u}_{\ep,t}\)=   u_{\ep} ^1 (y,u). \label{int.0ep}
\end{equation}
The local times $\wt L^{u}_{\ep, t} $ are jointly continuous in $u,t$.

   It follows from \cite[Lemma 2.1,
Theorems 4.2 and 5.2]{FMR} and the work of the last section  that
  for any closed interval $\De$ in $(0,1]$,
   \be
\lim_{h\to 0}\sup_{\stackrel{|u-v|\le h }{ u,v\in\De }} \frac{|\wt L_{\ep,t}^{u}-\wt L_{\ep,t}^{v}|}{  (|u-v| \log  1/| u-v|)^{1/2} }= \sup_{u\in \De}\(2\wt L_{\ep,t}^{u}\)^{1/2},  \quad a.e.\,\,\, t,  \,\,\, \wt P^{y}\,\,
a.s.\label{urev.m2}
\ee
for all $y>0$. 

In particular, if we set $I_{k,l}=[\frac{l-1}{2^{k}}, \frac{l+1}{2^{k}}]$, (overlapping intervals), 
   \bea
&&\lim_{h\to 0}\sup_{\stackrel{|u-v|\le h }{ u,v\in I_{k,l} }} \frac{|\wt L_{\ep,t}^{u}-\wt L_{\ep,t}^{v}|}{  (|u-v| \log  1/| u-v|)^{1/2} }= \sup_{u\in I_{k,l}}\(2\wt L_{\ep,t}^{u}\)^{1/2},\nn\\
&& \hspace{1 in}  \quad \forall k, 2\leq l\leq 2^{k}-1,\quad a.e.\,\,\, t,  \,\,\, \wt P^{y}\,\,
a.s.\label{urev.m2s}
\eea
for all $y>0$. 
Fix a  sequence $\ep_{n}\downarrow 0$. It follows from (\ref{urev.m2s}), the continuity of the local times and the density of $\mathcal{C}$ that
   \bea
&& \lim_{h\to 0}\sup_{\stackrel{|u-v|\le h }{ u,v\in [c_{j},1]\cap \mathcal{C}\cap I_{n,l} }} \frac{|\wt L_{ \ep_{n},t}^{u}-\wt L_{\ep_{n},t}^{v}|}{  (|u-v| \log  1/| u-v|)^{1/2} }= \sup_{u\in [c_{j},1]\cap \mathcal{C}\cap I_{n,l}}\(2\wt L_{\ep_{n},t}^{u}\)^{1/2},\nn\\
&& \hspace{1 in}\quad \forall n, 2\leq l\leq 2^{n}-1,\quad a.e.\,\,\, t,  \,\,\, \wt P^{y}\,\,
a.s.\label{urev.m3}
\eea
for all $y>0$.

We note in passing that it follows from \cite[p. 146]{Blu} that we can choose a sequence 
$\ep_{n}\downarrow 0$ such that for each fixed $t$, $ \wt   X_{\ep_{n},t}\to \wt   X_{t}$ with probability $1$.

 We note that $\wt L_{\ep,t}^{u}$ can only increase in $t$ when $\wt   X_{\ep,t}=u$, (or in the closure of such $t$'s). We consider two cases: when $t$ is in an excursion interval, that is  $t_{r}^{-}<t\leq t_{r}$ and when $t$ is in a holding interval, that is  $t_{r-1}<t\leq t_{r}^{-}$. (We note that $t_{r}^{-}, t_{r}$ are themselves functions of $\ep$).

   Fix $T<\ff$,   and let $W_{t}$ denote  standard Brownian motion on the real line. 
 It follows from 
  \cite[(14.5)]{book}, with $h=4,d=1,r=1$,
 \begin{equation}
 \|\sup_{\stackrel{ 0\leq s<t\leq T}{|t-s|\le 2^{-m}}}  |W_{t} -W_{s} | \,\, \|_{4}    \leq  C_{1}2^{-m/4},\label{3.klm1}
 \end{equation}
 so that for any $\al<1/4$
  \begin{equation}
 \|\sup_{\stackrel{ 0\leq s<t\leq T}{|t-s|\le 2^{-m}}}  \frac{ |W_{t} -W_{s} |}{2^{-\al m}}  \,\, \|_{4}    \leq C_{1}2^{-m(1/4-\al)},\label{3.klm2}
 \end{equation}
 and therefore
 \begin{equation}
 P\( \sup_{\stackrel{ 0\leq s<t\leq T}{|t-s|\le 2^{-m}}} \frac{ |W_{t} -W_{s} |}{2^{-\al m}}   \geq C_{3}     \)\leq C_{2} 2^{-m(1-4\al)}.\label{3.klm3}
 \end{equation}
 It follows that up to probability $C_{2} 2^{-m_{0}(1-4\al)}$ we have that 
 \begin{equation}
 \sup_{\stackrel{ 0\leq s<t\leq T}{|t-s|\le 2^{-m}}} \frac{ |W_{t} -W_{s} |}{2^{-\al m}}   \leq C_{3} , \qquad \forall m\geq m_{0}. \label{3.klm4}
 \end{equation}
 We will take $\al=1/8$, so that up to probability $C_{2} 2^{-m_{0}/2}$ we have that 
 \begin{equation}
 \sup_{\stackrel{ 0\leq s<t\leq T}{|t-s|\le 2^{-m}}} \frac{ |W_{t} -W_{s} |}{2^{- m/8}}   \leq C_{3} , \qquad \forall m\geq m_{0}. \label{3.klm5}
 \end{equation}
 We take $c_{j}=2^{-j}, \ep_{n}=2^{-8n}$ and $\de_{n}=\bb^{-1}_{\ep_{n}}=\rho_{n}2^{-\wt n}$ with $1\leq \rho_{n}\leq 2$.
 (The exact relationship of $n$ and $\wt n$ is not relevant for us. All we need is that $n\to \ff$ iff $\wt n \to \ff$.)
\medskip

  We now apply the above to differences of our process $\wt   X_{t}-\wt   X_{s},\, s<t\leq T$, which behave like the differences   $  W_{t}-W_{s},\, s<t\leq T$, of standard Brownian motion when  $\wt   X_{r}\in (0, \infty), s\leq r\leq t\leq T$.

(a) When  $t$ is in an excursion interval,   using (\ref{3.klm5}) we can choose $n$ large enough that if $X_{0}<c_{j+1}$ then also  $X_{ \ep_{n}}<c_{j}$ a.s.  so that by starting with  $\mu=\mu_{\ep_{n}}=:\bar\eta_{\ep_{n}}/\| \bar \eta_{\ep_{n}}   \|$, $X$ does not miss any points $u\in [c_{j},1]$. Thus in this case we have $\wt L_{ \ep_{n},t}^{u}=\wt L_{ 0,t}^{u}$ for all $u\in [c_{j},1]$.
\medskip

(a') On the other hand, if  $X_{0}=w\geq c_{j+1}$, (in particular if $w\geq c_{j}$), using  (\ref{3.klm5}) there will be an open interval $M_{ \ep_{n}, r}\subseteq (w-C_{3}\ep^{1/8}_{n},w+C_{3}\ep^{1/8}_{n})$ which may be missed by starting at $X_{ \ep_{n}}$.  
Thus for $t_{r}^{-}<t\leq t_{r}$ we have $\wt L_{ \ep_{n},t}^{u}=\wt L_{ 0,t}^{u}$ for all $u\in [c_{j},1]\cap M^{c}_{\ep_{n}, r}$.  Since $\ep^{1/8}_{n}=2^{-n}$, $M_{ \ep_{n}, r}$ will be contained in some fixed $C_{4}$ overlapping intervals of the form $I_{n,l}=[\frac{l-1}{2^{n}}, \frac{l+1}{2^{n}}]$.

How many excursion intervals will we have with  $X_{0}\geq c_{j+1}\geq \ep_{n}^{1/(8\cdot 16)}$ up to time $T$? Again using (\ref{3.klm5}), such an interval will have length at least $\ep_{n}^{1/16}$, hence up to time $T$ there will be at most $T\ep_{n}^{-1/16}$ such intervals. Thus, if we use $\wh M_{ \ep_{n}}$ to denote the union of the $T\ep_{n}^{-1/16}$ (random) open intervals $M_{\ep_{n}, r}\subseteq [c_{j},1]$, $\wh M_{ \ep_{n}}$ will be contained in $C_{4}T \ep_{n}^{-1/16}= C_{4}T2^{n/2}   $ overlapping intervals of the form $I_{n,l}=[\frac{l-1}{2^{n}}, \frac{l+1}{2^{n}}]$.
This will be a small proportion of the $  2^{n+1}$ intervals $I_{n,l}$, and this proportion $\to 0$ as $n\to \ff$.
\medskip

(b) When  $t$ is in a holding interval for $\wt   X_{\ep_{n},t}$, we have  $\wt   X_{\ep_{n},t}=0$.  The length of our holding interval is an exponential random variable with parameter $\bb_{\ep_{n}}=\| \bar  \eta_{\ep_{n}}   \|/p_{\ep_{n}}\uparrow \ff$. Set $\de_{n}=\bb^{-1}_{\ep_{n}}\downarrow 0$. The probability that a holding interval will have length $\geq \de^{1/4}_{n}$
is $e^{-1/\de^{3/4}_{n}}$. Let $\chi_{l}$ denote the length of the successive holding intervals. If  $V_{k}=\chi_{1}+\chi_{2}+\cdots +\chi_{k}$ then  $N_{T}=\sup\{k:\, V_{k}\leq T\}$  is Poisson with parameter $T/\de_{n}$. It follows that 
\begin{eqnarray}
&&P\(\text{up to time $T$ a holding interval will have length }\geq \de^{1/4}_{n}\)
\label{}
\\
&&=\sum_{k=0}^{\ff}P\(\text{up to time $T$ a holding interval will have length }\geq \de^{1/4}_{n}\,|\, N_{T}=k\)\nonumber\\
&&\hspace{4 in}\times P\( N_{T}=k\)
\nonumber\\
&&\leq \sum_{k=1}^{\ff}k   e^{-1/\de^{3/4}_{n}}   \frac{\(T/\de_{n}\)^{k}}{k!}e^{T/\de_{n}}
=\sum_{k=1}^{\ff}    e^{-1/\de^{3/4}_{n}}T/\de_{n}   \frac{\(T/\de_{n}\)^{(k-1)}}{(k-1)!}e^{T/\de_{n}}=Te^{-1/\de^{3/4}_{n}}/\de_{n}.\nonumber
\end{eqnarray}
This goes to $0$ as $\wt n\to \ff$.
\medskip

(b') Thus consider a holding interval of length $\leq \de^{1/4}_{n}$. What are we then missing about $\wt  X_{t}$?
 First of all, we are missing the first $\ep_{n}$ of each of our `jumping' intervals of length $>\ep_{n}$, but this has been dealt with in (a) and (a'). Secondly we are missing  excursion intervals of $\wt X$ of length $\leq \de^{1/4}_{n}$. But since they have to end at zero, by (\ref{3.klm5}), for $n$ sufficiently large, we will never have $\wt X_{t}\geq c_{j}$.

As mentioned, if for some $u,t, n$ we have that   $\wt L_{ \ep_{n'},t}^{u}=\wt L_{ \ep_{n},t}^{u}$ for all $ n'\geq n$, we set $\wt L_{ 0,t}^{u}=\wt L_{ \ep_{n},t}^{u}$. We will also write this as $\wt L_{ \ep_{n}, \ast,t}^{u}$ to emphasize that this doesn't change as $n$ increases. For which $u,t$ 
does $\wt L_{ 0,t}^{u}$ exist? It is clear from the above that up to some probability that goes to $0$,  the only problem will come from (a'). The only $u,t$ for which 
 $\wt L_{ 0,t}^{u}$ may not exist  are the $u,t $ which are the initial points of the   excursions of $\wt X$. Note that such $u\notin \mathcal{C}$ by (\ref{lil.c0}). In particular, we will see that $\wt L_{0,t}^{u}$ will be defined for each $u\in \mathcal{C}$. In more detail:
 
 It follows from (\ref{urev.m3})  and what we have just shown that the event
 \begin{equation}
\lim_{h\to 0}\sup_{\stackrel{|u-v|\le h }{ u,v\in (\wh M_{ \ep_{n}})^{c}\cap[c_{j},1]\cap \mathcal{C} }} \frac{|\wt L_{ \ep_{n}, \ast,t}^{u}-\wt L_{ \ep_{n}, \ast,t}^{v}|}{  (|u-v| \log  1/| u-v|)^{1/2} }= \sup_{u\in (\wh M_{ \ep_{n}})^{c}\cap [c_{j},1]\cap \mathcal{C}}\(2\wt L_{ \ep_{n}, \ast,t}^{u}\)^{1/2},   \label{}
 \end{equation}
 for a.e. $t\leq T$,
 has $\wt P^{y}$ probability $=1-o_{n}(1)$, for all $y>0$. 
Taking the limit as $ n\to \ff$ and then  $T\to \ff$     we obtain 
   \be
\wt P^{y}\(\lim_{h\to 0}\sup_{\stackrel{|u-v|\le h }{ u,v\in [c_{j},1]\cap \mathcal{C} }} \frac{|\wt L_{0,t}^{u}-\wt L_{0,t}^{v}|}{  (|u-v| \log  1/| u-v|)^{1/2} }=  \sup_{u\in  [c_{j},1]\cap \mathcal{C}}\(2\wt L_{0,t}^{u}\)^{1/2},  \quad a.e.\,\,  t \)= 1 .\label{lil.44a}
\ee

The analysis above shows that that for each $u\in \mathcal{C}$, $\wt L_{ 0,t}^{u}$ is a CAF in $t$ and in particular a local time for $u$.

Next,  recall that our local times $\wt L_{ \ep_{n},t}^{u}$ are normalized so that for any   $x, u>0$,
 \begin{equation}
  \wt E^{ x}\(\int_{0}^{\ff}  e^{-  s}\, d_{s}\wt L^{u}_{ \ep_{n},s}\)=   u_{\ep_{n}}^1 (x,u). \label{maz.1vw}
\end{equation}
It follows from the above and  (\ref{int.0dense})   that for all $x>0$ and $u\in \mathcal{C} $
 \bea
  \wt E^{ x}\( \int_{0}^{\ff}  e^{- s}\, d_{s}\wt L^{u}_{0,s}\)&=& \lim_{n\to\ff}  u_{\ep_{n}}^1 (x,u)=u^1 (x,u),\label{maz.1vwa}
\eea
which is
 (\ref{int.0ep0}) and this completes the proof of Theorem \ref{theo-uniformmodulus}   in the  case when $\lim_{\ep\to 0}\|\bar\eta_{\ep}\|=\ff$. Sub-intervals of $(0,1]$ can be treated similarly.

 We next consider the  case when $\lim_{\ep\to 0}\|\bar\eta_{\ep}\|<\ff.$ This requires that $q=0$ and $\| \eta \|<\ff.$ It then follows from \cite[Chapter II, Theorem 3.8]{Blu} that $p>0$.
 
 We can then take $\eta_{\ep}=\eta$ so that with $\la=1$, (\ref{blu.17}) becomes 
   \begin{equation}
U_{\bullet}^{1}h(x)=V^{1}h(x) +e^{-\sqrt{2}\,\, x}\,\,\frac{ph(0)+  \eta V^{1}h}{p+\int_{0}^{\ff}\(1-e^{-\sqrt{2 }\,\, x}\)\,d\eta (x)}.   \label{blu.17p}
 \end{equation}
 (Since this does not depend on $\ep$, we use the subscript $\bullet$, but in fact this is the correct resolvent. We don't have to work with approximations!)
 
 In this situation we have $p+\int_{0}^{\ff}\(1-e^{-\sqrt{2 }\,\, x}\)\,d\eta (x)=1$, see \cite{Blu}, so that we have 
    \begin{equation}
U_{\bullet}^{1}h(x)=V^{1}h(x) +e^{-\sqrt{2}\,\, x}\,\,\( ph(0)+  \eta V^{1}h\).   \label{blu.17p1}
 \end{equation}
 In particular, for all $x,y>0$
     \begin{equation}
u_{\bullet}^{1}(x,y)=v^{1}(x,y)+e^{-\sqrt{2}\,\, x}\,\,  \int v^{1}(z,y)\,  \eta (dz).   \label{blu.17p2}
 \end{equation}
 The precise uniform modulus of continuity for the local times $\wt L^{u}_{s}$, $u\in (0,1]$ follows immediately from \cite[Lemma 2.1,
Theorems 4.2 and 5.2]{FMR}  and the work of this section.

We note that in our work which uses Isomorphism Theorems to derive uniform moduli for local times  we restrict to compact sets $K$, see \cite[section 9.5]{book} and \cite{FMR}. This is because we need to use the result that for any continuous Gaussian process $G_{x}$ and any  $\ep>0$, the set where $\sup_{x\in K}|G_{x}(\om)|\leq \ep$ has positive probability. This uses \cite[Lemma 5.3.5]{book}. However, in our case, the
continuous Gaussian processes $\eta_{i, 0} (x)$ with covariance $ v^{0}(x,y)$ and 
 $\eta_{p} (x)$  with covariance $  v^{p}(x,y)$ are $0$ when $x=0$, so that we can obtain  uniform moduli for local times on $(0,1]$.

Before closing this section we want to collect some facts which will be useful for Section \ref{sec-lt0}. By \cite[Theorem 3]{R83} with $p_{1}=0$, $p_{2}=q/\sqrt{2}$, $p_{3}=p$ and $p_{4}=\eta$, we have 
\begin{equation}
U^{\la}h(x)=V^{\la}h(x) +e^{-\sqrt{2\la}\,\, x}\,\,\frac{ph(0)+ q\sqrt{2}\int_{0}^{\ff} e^{-\sqrt{2\la}\,\, y}h(y)\,dy +\eta V^{\la}h}{p\la+q\sqrt{\la}+\int_{0}^{\ff}\(1-e^{-\sqrt{2\la}\,\, x}\)\,d\eta(x)},\label{blu.42r}
\end{equation}
see  (\ref{r.42r}).   Since  $\sqrt{\la}\, R^{\la}h(0)=\sqrt{2}\int_{0}^{\ff} e^{-\sqrt{2\la}\,\, y}h(y)\,dy$ we can write this as 
\begin{equation}
U^{\la}h(x)=V^{\la}h(x) +e^{-\sqrt{2\la}\,\, x}\,\,\frac{ph(0)+ q\sqrt{\la}\, R^{\la}h(0) +\eta V^{\la}h}{p\la+q\sqrt{\la}+\int_{0}^{\ff}\(1-e^{-\sqrt{2\la}\,\, x}\)\,d\eta(x)}.\label{blu.42}
\end{equation}
When $p=0$ this becomes
\begin{equation}
U^{\la}h(x)=V^{\la}h(x) +e^{-\sqrt{2\la}\,\, x}\,\,\frac{ q\sqrt{\la}\, R^{\la}h(0) +\eta V^{\la}h}{ q\sqrt{\la}+\la \eta V^{\la}1}.\label{blu.42b}
\end{equation}
 We note that if, in addition, we take $ \eta=0$, this becomes
 \begin{equation}
U^{\la}h(x)=V^{\la}h(x) +e^{-\sqrt{2\la}\,\, x}\,\,\, R^{\la}h(0)=R^{\la}h(x),\label{blu.42c}
\end{equation}
 as should be. 
 
Returning to (\ref{blu.42r}) but with $p=0$, it follows that $U^{\la} $ has a density 
 \begin{equation}
u^{\la}(y,u)=    v^{\la}(y,u)  +e^{-\sqrt{2\la}\,\, y}\,\,\frac{ q\sqrt{2} e^{-\sqrt{2\la}\,\, u}+\int_{0}^{\ff}v^{\la}(z,u) \,d\eta(z)}{ q\sqrt{\la}+\int_{0}^{\ff}\(1-e^{-\sqrt{2\la}\,\, z}\)\,d\eta(z)}\label{blu.42dens}
\end{equation}
 with respect to Lebesgue measure. In particular
  \begin{equation}
u^{\la}(y,0)=     e^{-\sqrt{2\la}\,\, y}\,\,\frac{ q\sqrt{2} }{ q\sqrt{\la}+\int_{0}^{\ff}\(1-e^{-\sqrt{2\la}\,\, z}\)\,d\eta(z)}.\label{blu.42dens0}
\end{equation}

 \section{Exact uniform moduli of continuity for the local times of `Feller Diffusions`}\label{sec-FD}

Let  $\mathcal{Z}$ be a  diffusion in $R^1$ that is regular and without traps  and  is symmetric with respect to a $\si$--finite measure  $m$, called the speed measure, which is absolutely continuous with respect to Lebesgue measure.  
(A diffusion is   regular and without traps when
$P^{x}\( T_{y}<\ff\)>0,   \forall x,y\in  R^1$.) We   consider diffusions $\mathcal{Z}$ with  generators of the form
\begin{equation}
\wt  L=\frac{1}{2}a^{2}(x)\frac{d^{2}}{dx^{2}}+c(x)\frac{d}{dx}\label{gendef.1},
\end{equation}
where $a^{2}(x)$ and $c(x)\in C\(R^1\)$ and $a(x)\ne 0$ for any $x\in R^1$.    For details see   \cite[Section 7.3]{RY} and \cite[Chapter 16]{Breiman}.

   Let  $s(x)$ be a scale function for $\cal Z$.  It is   strictly  increasing and  unique up to a linear transformation.  Since $a^{2}(x)$ and $c(x)\in C\(R^1\)$, $s(x)\in C^2(R^1)$.   We take $s(0)=0$  and consider the case where $\lim_{x\to \ff}s(x)=\ff$.

 For $\la>0$, the $\la-$potential density of $\cal Z$ with respect to  the speed measure     is,  
 \begin{equation} \label{diff.1}
\ov u^\la (x,y)= \left\{
 \begin{array} {cc}
 p_\la (x)q_\la (y),& \quad x\leq y  
 \\
 q_\la(x)p_\la(y),& \quad y\leq x   
\end{array}  \right. ,
\end{equation}
where $p_\la$ and $q_\la$ are in $ C^{2}(R^{1})$ and  are positive and    $p_\la$ is strictly  increasing and  
$q_\la$ is strictly  decreasing.     See   \cite[(4.114)]{book}.  

We have
\begin{equation}
E^{x}\(e^{-\la \si}\)=\frac{\ov u^{\la}(x,0 )}{\ov u^{\la}(0,0)},\label{diff.2}
\end{equation}
and for all $\la>0$,
\begin{equation} \label{5.44}
   s'(x)=q_\la (x)p'_\la (x)-p_\la (x)q'_\la (x).
\end{equation}
 See e.g., \cite[Theorem 9.1]{MRLIL}.
 
Let $(\wt X_{t}, \wt P^{x})$ be a recurrent Borel right process  with state space $[0,\ff)$ which behaves like $\mathcal{Z}$ up until the hitting time $\si$ of the point $0$. $\wt X$ has a $1$-potential operator $U^{1}$ which, restricted to $(0,\ff)$, has a continuous density $u ^1 (x,u)$ with respect to  speed measure. Let $\wt L^{u}_{s}$ denote the local times of $\wt X$, normalized so that for
 $y, u>0$,
 \begin{equation}
  \wt E^{ y}\(\int_{0}^{\ff}  e^{-  s}\, d_{s}\wt L^{u}_{s}\)=   u ^1 (y,u). \label{diff.0}
\end{equation}

 The goal of this section is to prove the following Theorem.
 
 \bt\label{theo-diffuniformmodulus} For any $y>0$
   \be
\lim_{h\to 0}\sup_{\stackrel{|u-v|\le h }{ u,v\in (0,1] }} \frac{|\wt L_{t}^{u}-\wt L_{t}^{v}|}{  (|s(u)-s(v)| \log  1/|u-v |)^{1/2} }= \sup_{u\in (0,1]}\(2\wt L_{t}^{u}\)^{1/2},  \quad a.e.\,\,\, t,  \,\,\, \wt P^{y}\,\,
a.s.\label{diff.00}
\ee
 \et

 Let  $X=\{X_{t},t\ge 0 \}$, be the process with state space  $S=(0,\ff)$ that is obtained by starting  $\mathcal{Z}$ in $S$  and then  killing it  the first time it hits 0. The $\la$- potential density of $X$ with respect to the speed measure  $m$ for $\la>0$ is,
\begin{equation} \label{d5.6}
v^{\la}(x,y)=\ov u^{\la}(x,y)-\frac{\ov u^{\la}(x,0 )\ov u^{\la}(0, y)}{\ov u^{\la}(0,0)}=\ov u^{\la}(x,y)-E^{x}\(e^{-\la \si}\)\ov u^{\la}(0, y);
\end{equation} 
see  \cite[4.165]{book},
and 
  \begin{equation} \label{diff.3}
v^{0}(x,y)=s(x)\wedge s(y),\qquad  x,y>0.
\end{equation}
(We get this by taking the limit $a $ goes to 0 and  $b $ goes to infinity in \cite[VII. Corollary 3.8]{RY}.)

We note that reflecting our diffusion motion in $[0,\ff)$ has continuous strictly positive  $\la-$potential   densities   with respect to $m$ on $[0,\ff)$:
\be
 r^{\la}(x,y)= \ov u^{\la}(x,y)+\frac{\ov u^{\la}(x,0 )\ov u^{\la}(0, y)}{\ov u^{\la}(0,0)}, \quad x,y\geq 0. \label{diff.4}
\ee 
For measurable functions  $h$ on $[0,\ff)$ we  define  the $\la$--th potential operator,
\begin{equation}
R^{\la}h(x)=\int_{0}^{\ff}   r^{\la}(x,y)h(y) \,dm (y).    \label{diff.5}
\end{equation}
In particular
\begin{equation}
R^{\la}h(0)=\int_{0}^{\ff}   r^{\la}(0,y)h(y) \,dm (y)=2\int_{0}^{\ff}  \ov u^{\la}(0,y)h(y) \,dm (y).    \label{diff.6}
\end{equation}
It then follows that 
\begin{equation}
V^{\la}h(x)+ E^{x}\(e^{-\la \si}\)  R^{\la}h(0)=R^{\la}h(x).\label{diff.7}
\end{equation}

\bl\label{lem-difref}
\begin{equation}
\lim_{x\to 0}\frac{V^{\la}h(x)}{s(x)}=\frac{1}{2 \ov u^{\la}(0,0) } \,R^{\la}h(0).\label{diff.8}
\end{equation}
\el

Proof: 
\begin{eqnarray}
&&V^{\la}h(x)=\int _{0}^{\ff}\ov u^{\la}(x,y)h(y)\,dm(y)-\frac{\ov u^{\la}(x,0 ) }{\ov u^{\la}(0,0)}\int _{0}^{\ff}\ov u^{\la}(0,y)h(y)\,dm(y)
\nn
\\
&&=\int _{0}^{x}\ov u^{\la}(x,y)h(y)\,dm(y)+\int _{x}^{\ff}\ov u^{\la}(x,y)h(y)\,dm(y)
\label{diff.9}\\
&&-\frac{\ov u^{\la}(x,0 ) }{\ov u^{\la}(0,0)}\int _{0}^{x}\ov u^{\la}(0,y)h(y)\,dm(y)-\frac{\ov u^{\la}(x,0 ) }{\ov u^{\la}(0,0)}\int _{x}^{\ff}\ov u^{\la}(0,y)h(y)\,dm(y)\nn\\
&&=q_\la (x)\int _{0}^{x} p_\la (y)h(y)\,dm(y)+p_\la (x)\int _{x}^{\ff} q_\la (y)h(y)\,dm(y)
\nonumber\\
&&-\frac{q_\la (x) }{q_\la (0)}\int _{0}^{x}\ov u^{\la}(0,y)h(y)\,dm(y)-\frac{q_\la (x) }{q_\la (0)}\int _{x}^{\ff}\ov u^{\la}(0,y)h(y)\,dm(y).\nn
\end{eqnarray}

Consider first
\begin{eqnarray}
&&I=p_\la (x)\int _{x}^{\ff} q_\la (y)h(y)\,dm(y)-\frac{q_\la (x) }{q_\la (0)}\int _{x}^{\ff}\ov u^{\la}(0,y)h(y)\,dm(y)
\nn
\\
&&=\(\frac{p_\la (x) }{p_\la (0)}-\frac{q_\la (x) }{q_\la (0)}\)\int _{x}^{\ff}\ov u^{\la}(0,y)h(y)\,dm(y).
\label{diff.10}
\end{eqnarray}
Using the fact that $p_\la$ and $q_\la$ are in $ C^{2}(R^{1})$, for small $x$ we can write $p_\la (x)=p_\la (0)+x p'_\la (0)+O(x^{2}) $ and $q_\la (x)=q_\la (0)+x q'_\la (0)+O(x^{2})$ so that by (\ref{5.44})
\begin{eqnarray}
&&
\(\frac{p_\la (x) }{p_\la (0)}-\frac{q_\la (x) }{q_\la (0)}\)=(x+O(x^{2}))\(\frac{p'_\la (0) }{p_\la (0)}-\frac{q'_\la (0) }{q_\la (0)}\)\nn\\
&&=(x+O(x^{2})) \frac{p'_\la (0)q_\la (0)- q'_\la (0)p_\la (0)}{p_\la (0)q_\la (0)}\sim \frac{xs'(0)}{\ov u^{\la}(0,0)}\sim \frac{s(x)}{\ov u^{\la}(0,0)}, \label{diff.11}
\nonumber
\end{eqnarray}
since $s(0)=0$.
Hence, using (\ref{diff.6})
\begin{equation}
\lim_{x\to 0}\frac{I}{s(x)}=\frac{1}{2 \ov u^{\la}(0,0) } \,R^{\la}h(0).\label{diff.12}
\end{equation}

Now consider 
\begin{eqnarray}
&&II=q_\la (x)\int _{0}^{x} p_\la (y)h(y)\,dm(y)-\frac{q_\la (x) }{q_\la (0)}\int _{0}^{x}\ov u^{\la}(0,y)h(y)\,dm(y)
\label{diff.13}
\\
&&=q_\la (x)\(\int _{0}^{x} p_\la (y)h(y)\,dm(y)-\frac{1 }{q_\la (0)}\int _{0}^{x}  p_\la (0) q_\la (y)h(y)\,dm(y)\)
\nonumber\\
&&=\frac{q_\la (x) }{q_\la (0)}\(\int _{0}^{x} \(p_\la (y)q_\la (0)-  p_\la (0) q_\la (y)\)h(y)\,dm(y)\).
\nonumber 
\end{eqnarray}
As before, using $  \(p_\la (y)q_\la (0)-  p_\la (0) q_\la (y)\)\sim  s(y)$ for small $y\leq x$ we see that 
\begin{equation}
\lim_{x\to 0}\frac{II}{s(x)}=0.\label{diff.14}
\end{equation}
\qed

As in the case of Feller Brownian motions we will approximate our local times 
by the local times of the holding and jumping process described in (\ref{blu.1}) with $\la$-potential operator
  \begin{equation}
U_{\ep}^{\la}h(x)=: V^{\la}h(x) +E^{x}\( e^{-\la \si }\)\frac{p_{\ep} h(0) +  \bar \eta_{\ep} V^{\la}h} {   \la p_{\ep}+ \la  \bar \eta_{\ep} V^{\la}1  }.   \label{diff.15}
 \end{equation}

 
 In  (\ref{diff.15}) take  
 \begin{equation}
\bar\eta_{\ep}=\frac{q}{s(\ep)}\,\de_{\ep}+\eta_{\ep}\label{diff.15a}
\end{equation}
with $\eta_{\ep}$ as in (\ref{elaw.1}) to obtain  
  \begin{equation}
U_{\ep}^{\la}h(x)=V^{\la}h(x) +E^{x}\( e^{-\la \si }\)\,\,\frac{p_{\ep}h(0)+ \frac{q}{s(\ep)}\, V^{\la}h(\ep) + \eta_{\ep} V^{\la}h}{\la p_{\ep}+ \frac{q \la }{s(\ep)}\, V^{\la}1(\ep) + \la \eta_{\ep} V^{\la}1}.   \label{diff.17q}
 \end{equation}
Considering first the case where $p_{\ep}=p>0$ and using Lemma \ref{lem-difref}
we find that
  \begin{equation}
\lim_{\ep\to 0}U_{\ep}^{\la}h(x)=V^{\la}h(x) +E^{x}\( e^{-\la \si }\)\,\,\frac{p h(0)+ \frac{q}{2 \ov u^{\la}(0,0) } \,R^{\la}h(0) + \eta  V^{\la}h}{\la p + \frac{q \la }{2 \ov u^{\la}(0,0) } \,R^{\la}1(0) + \la \eta  V^{\la}1}.   \label{diff.18}
 \end{equation}
In particular with $\la=1$ we have
  \begin{equation}
\lim_{\ep\to 0}U_{\ep}^{1}h(x)=V^{1}h(x) +E^{x}\( e^{-  \si }\)\,\,\frac{p h(0)+ \frac{q}{2 \ov u^{1}(0,0) } \,R^{1}h(0) + \eta  V^{1}h}{ p + \frac{q}{2 \ov u^{1}(0,0) } \,R^{1}1(0) +\eta  V^{1}1}.   \label{diff.19}
 \end{equation}
When $p_{\ep}\to 0$ we obtain a similar expression with $p=0.$

As before, let $\wt L_{\ep,t}^{u}$ denote the local times of the process  
 $\wt   X_{\ep,t}=\wt X(v^{p}, m ,\mu_{\ep},\bb_{\ep})$.    
   It follows from \cite[Lemma 2.1,
Theorems 4.2 and 5.4]{FMR} and the work of the Section 2  that
  for any closed interval $\De$ in $(0,1]$,
   \be
\lim_{h\to 0}\sup_{\stackrel{|u-v|\le h }{ u,v\in\De }} \frac{|\wt L_{\ep,t}^{u}-\wt L_{\ep,t}^{v}|}{  (|s(u)-s(v)| \log  1/| u-v|)^{1/2} }= \sup_{u\in \De}\(2\wt L_{\ep,t}^{u}\)^{1/2},  \quad a.e.\,\,\, t,  \,\,\, \wt P^{y}\,\,
a.s.\label{diff.20}
\ee
for all $y>0$. 
As before, using [1, Proposition 6.11] to obtain the analogue of (3.31), we
can use this to prove Theorem 4.1.

 \section{Joint continuity of the local time on $\mathbf{[0,1]\times R_{+}^{1}}$}\label{sec-lt0}
 
 In Feller Brownian motion the local time at $0$ plays a special role. The reader will have noticed that Theorem \ref{theo-uniformmodulus}, our exact uniform moduli of continuity result,  is for $\wt L^{u}_{t}$ with $u\in (0,1]$. When $p\neq 0$, so that $0$ is `sticky', the 1-potential density of $\wt X$ with respect to Lebesgue measure is not continuous at $0$. When $p=0$ but $ q>0$, that is, when the Feller Brownian motion has some reflecting  Brownian motion, we will give examples below where $\wt L^{u}_{t}$ is jointly continuous on $[0,1]\times R_{+}^{1}$. The case where $p=0, q=0$ is the really interesting case. This is because the 1-potential density of $\wt X$ with respect to Lebesgue measure, $u^1(y,u)$, is $0$ when $u=0$. We cannot normalize $\wt L^{0}_{t}$ by $u^1(y,u)$. We have to normalize differently, and then we can give examples where the local time is jointly continuous on $[0,1]\times R_{+}^{1}$.
 
 \subsection{$\mathbf{p=0}$ but $\mathbf{ q>0}$}
 
  When $ \eta=0$, by (\ref{blu.42c}) we have reflecting Brownian motion which has a jointly continuous local time $\wt L^{u}_{t}$ normalized by 
   \begin{equation}
  \wt E^{ y}\(\int_{0}^{\ff}  e^{-  t}\, d_{t}\wt L^{u}_{t}\)=   r ^1 (y,u),\qquad y,u\in [0,\ff). \label{3.50}
\end{equation}

 We next consider the case where $p=0$, $ q>0$  and $ \eta\neq 0$.  It suffices to consider 
 the local times $\wt L^{u}_{t}$ for $u\in [0,1]$, normalized by 
    \bea
  &&\wt E^{ y}\(\int_{0}^{\ff}  e^{- t}\, d_{t}\wt L^{u}_{t}\)=   u ^1 (y,u) \label{3.51}\\
  &&\hspace{1 in}=
v^1 (y,u) +e^{-\sqrt{2}\,\, y}\(q \sqrt{2}e^{-\sqrt{2}\,\, u}+\,\,\int_{0}^{\ff}v^1 (z,u)\,d\eta(z)\),\nn
\eea
with
\begin{equation}
q+\int_{0}^{\ff}\(1-e^{-\sqrt{2}\,\, z}\)\,d\eta(z)=1,\label{3.51a}
\end{equation}
see (\ref{blu.42r}).
We will provide conditions on $\eta$ which will guarantee that $\wt L^{u}_{t}$  is jointly continuous on $[0,1]\times R_{+}^{1}$.  It will then follow that  Theorem \ref{theo-uniformmodulus} extends to $[0,1]$. To see this write (\ref{int.1}) as
   \be
\lim_{h\to 0}\sup_{\stackrel{0<|u-v|\le h }{ u,v\in (0,1] }} \frac{|\wt L_{t}^{u}-\wt L_{t}^{v}|}{  (|u-v| \log  1/| u-v|)^{1/2} }= \sup_{u\in (0,1]}\(2\wt L_{t}^{u}\)^{1/2},  \quad a.e.\,\,\, t,  \,\,\, \wt P^{y}\,\,
a.s.\label{com.1}
\ee

Then by fixing $v>0$ and letting $u\to 0$, we see that for each $h>0$ we have 
\begin{equation}
\sup_{\stackrel{0<|u-v|\le h }{ u,v\in (0,1] }} \frac{|\wt L_{t}^{u}-\wt L_{t}^{v}|}{  (|u-v| \log  1/| u-v|)^{1/2} }=\sup_{\stackrel{0<|u-v|\le h }{ u,v\in [0,1] }} \frac{|\wt L_{t}^{u}-\wt L_{t}^{v}|}{  (|u-v| \log  1/| u-v|)^{1/2} },\label{com.2}
\end{equation}
and of course $\sup_{u\in (0,1]}\(2\wt L_{t}^{u}\)^{1/2}=\sup_{u\in [0,1]}\(2\wt L_{t}^{u}\)^{1/2}$. 

We now proceed to describe our conditions.
Note that with $q>0$,   $ u ^1 (y,u)$ is strictly positive for $y,u\geq 0$. In fact 
\begin{equation}
u ^1 (y,u)\geq   e^{-\sqrt{2}\,\, y}q \sqrt{2}e^{-\sqrt{2}\,\, u}.\label{3sp}
\end{equation}
(\ref{3.51a}) implies that 
\begin{equation}
\int_{0}^{\ff}z\wedge 1\,d\eta(z)<\ff.\label{3.51aa}
\end{equation}
Note that for $u\leq y\leq 1$
\begin{equation}
v^{1}(y,u)=e^{-\sqrt{2}\,\,y}\(\frac{e^{\sqrt{2}\,\,u}-e^{-\sqrt{2}\,\,u}}{\sqrt{2}}\)\leq cu=c(y \wedge u). 
\label{3.92}
\end{equation}
and it follows from this that 
\begin{equation}
v^{1}(y,u)\leq c(y \wedge u \wedge 1).\label{3base}
\end{equation}
It then follows  that
\begin{equation}
g(u)=:q \sqrt{2}e^{-\sqrt{2}\,\, u}+\,\,\int_{0}^{\ff}v^1 (z,u)\,d\eta(z)\label{3.53}
\end{equation}
is bounded by a constant independent of $u$. 

Let $f(x)=e^{-\sqrt{2}\,\, x}$ and $\Ga (x,y)=v^1(x,y)+f(x)g(y)$. Note that
\begin{equation}
\sup_{x,y \in [0,1]}\Ga (x,y)<\ff\label{3.900a}
\end{equation}
and
\begin{equation}
\inf_{x,y \in [0,1]}\Ga (x,y)>0.\label{3.900b}
\end{equation}
Let
\begin{equation}
d^{2}(x,y)=\Ga^{2} (x,x)+\Ga^{2} (y,y)-2\Ga (x,y)\Ga (y,x),\label{3.90}
\end{equation}
and for any probability measure $\mu$ let
\begin{equation}
J_{T, d,\mu}(a)=\sup_{t\in T}\int_{0}^{a} \(\log \frac {1}{\mu\(B_{d}(t,u)\)}\)^{1/2}\,du,\label{3.90A}
\end{equation}
 where $B_d(t,u)=\{x\mid d(x,t)\le u\}$ is the $d$-ball of radius $u$ centered at $t$.

It follows from \cite[Theorem 1.2, Lemma 5.1, (73), (74)]{MRsuf} and \cite[Lemma 4.2]{KMR}, that $d(x,y)$ is a metric  and that $\wt L^{u}_{t}$ will be  jointly continuous on $[0,1]\times R_{+}^{1}$ if for some probability measure $\mu$
\begin{equation}
\lim_{\ep\to 0}J_{[0,1], d,\mu}(\ep)=0.\label{3.90a}
\end{equation}
(We point out that \cite[Theorem 1.2]{MRsuf} is meant to obtain joint continuity over 
$R_{+}^{1}\times R_{+}^{1}$ so that it requires (\ref{3.90a}) for $J_{K, d,\mu}$ for all compact $K\subseteq R_{+}^{1}$, but the proof shows that if we only have (\ref{3.90a})
for some compact $K$ we get joint continuity over 
$K\times R_{+}^{1}$).

We have
\begin{eqnarray}
&&\Ga^{2} (x,x)+\Ga^{2} (y,y)-2\Ga (x,y)\Ga (y,x)
\label{3.64}
\\
&&=\Ga  (x,x) \(\Ga  (x,x)-\Ga (y,x)\)+\Ga (y,x) \(\Ga  (x,x)-\Ga (x,y)\)
\nonumber\\
&&+\Ga  (y,y) \(\Ga  (y,y)-\Ga (y,x)\)+\Ga (y,x) \(\Ga  (y,y)-\Ga (x,y)\),
\nonumber
\end{eqnarray}
so   to bound (\ref{3.64}) we only need to bound  expressions such as 
\begin{eqnarray}
&&\Ga  (x,x)-\Ga (y,x)
\label{3.65}
\\
&&=\(v^1  (x,x)-v^1  (y,x)\)+(f(x)-f(y))g(x)
\nonumber
\end{eqnarray}
and 
\begin{eqnarray}
&&\Ga  (x,x)-\Ga (x,y)
\label{3.66}
\\
&&=\(v^1   (x,x)-v^1  (x,y)\)+f(x)(g(x)-g(y)).
\nonumber
\end{eqnarray}

We provide an increasingly more general sequence of conditions, A, B and C.  In our applications of (\ref{3.90a}) we take the probability measure $\mu$
to be a multiple of Lebesgue measure. 
\medskip

{\bf Condition A:}

Assume that for some $\de>0$
\begin{equation}
\int_{0}^{1}z^{1-\de}\,d\eta(z)<\ff,\label{3.51aad}
\end{equation}
so that 
\begin{equation}
\int_{0}^{\ff}z^{1-\de}\wedge 1\,d\eta(z)=C_{\de}<\ff.\label{3.501}
\end{equation}
Then, using the fact that 
\begin{equation}
|v^1 (z,x)-v^1 (z,y)|\leq 2|x-y|,\label{3.501a}
\end{equation}
 and (\ref{3base}) we see that 
\begin{eqnarray}
&&|    \int_{0}^{\ff}\(v^1 (z,x)-v^1 (z,y)\)\,d\eta(z)      |
\label{3.501b}
\\  
&& \leq  \int_{0}^{\ff}|v^1 (z,x)-v^1 (z,y)|\,d\eta(z)
\nonumber\\  
&& \leq  \int_{0}^{\ff}|v^1 (z,x)-v^1 (z,y)|^{\de}|v^1 (z,x)-v^1 (z,y)|^{1-\de}\,d\eta(z)
\nonumber\\  
&& \leq  4|x-y|^{\de}\int_{0}^{\ff}z^{1-\de}\wedge 1\,d\eta(z)=4C_{\de}|x-y|^{\de},
\nonumber
\end{eqnarray}
and clearly
\begin{equation}
|f\(x\)-f\(y\)|\leq c|x-y|,\label{}
\end{equation}
and we have (\ref{3.501a}).
Hence 
\begin{equation}
|g(x)-g(y)|\leq C |x-y|^{\de}. \label{}
\end{equation}
It follows that $d(x,y)\leq   C |x-y|^{\de/2}$, and 
\bea
B_{d}(t,u)& = &\{x\,|\,d(x,t)\leq u\}\supseteq \{x\,|\,  |x-t|^{\de/2}\leq cu\}\label{3.58lf}\\
&=&   \{x\,|\, |x-t|\leq c' u^{2/\de}\}, \nn
\eea
so that if $\mu$
is a multiple of Lebesgue measure,  $\mu (B_{d}(t,u))\geq  c''u^{2/\de}$ and therefore
\begin{equation}
\log \frac{1}{\mu (B_{d}(t,u))}\leq c'\log u.\label{3.59lf}
\end{equation}
Clearly (\ref{3.90a}) holds, so that $\wt L^{u}_{t}$ is jointly continuous on $[0,1]\times R_{+}^{1}$. 
\medskip

{\bf Condition B:}

More generally, assume that
 for some $\ga>1$
\begin{equation}
\int_{0}^{1}z\,\log ^{\ga}(1/z)\,d\eta(z)<\ff,\label{3.51aal}
\end{equation}
so that with $h(t)=t\,\log_{\ast} ^{\ga}(1/t)$, where $\log_{\ast}(x)=\log(x)$ if $x\geq 1$ and 
$\log_{\ast}(x)=1$ if $x<1 $, 
\begin{equation}
\int_{0}^{\ff}h\(z\)\wedge 1\,d\eta(z)=c_{h}<\ff.\label{3.501l}
\end{equation}
It then follows as before that 
\begin{equation}
|g(x)-g(y)|\leq c \log_{\ast}  ^{-\ga}(1/|x-y|).  \label{3.57l}
\end{equation}
It follows that $d(x,y)\leq c \log_{\ast}  ^{-\ga/2}(1/|x-y|)$, and 
then in the notation of \cite[Theorem 1.2]{MRsuf}
\bea
B_{d}(t,u)&= &\{x\,|\, d(x,t)\leq u\}\supseteq \{x\,|\, \log_{\ast}  ^{-\ga/2}(1/|x-t|)\leq cu\}\label{3.58l}\\
&=&   \{x\,|\, \log_{\ast}  ^{\ga/2}(1/|x-t|)\geq \frac{1}{cu}\}=  \{x\,|\, \log_{\ast}(1/|x-t|)\geq c' u^{-2/\ga}\} \nn\\
&=&   \{x\,|\,  1/|x-t|\geq e^{c' u^{-2/\ga} }\}=  \{x\,|\, |x-t|\leq e^{-c'u^{-2/\ga} }\}, \nn
\eea
so that $\mu (B_{d}(t,u))\geq e^{-c'u^{-2/\ga} }$ and therefore
\begin{equation}
\log \frac{1}{\mu (B_{d}(t,u))}\leq c'u^{-2/\ga}.\label{3.59l}
\end{equation}
By (\ref{3.90a}) we will have continuity if
\begin{equation}
\lim_{\ep\to 0}\int_{0}^{\ep} \{c'u^{-2/\ga}\}^{1/2}\,du=c''\lim_{\ep\to 0}\int_{0}^{\ep}  \frac{1}{u^{1/\ga}} \,du=0,\label{3.59r}
\end{equation}
which will happen with    $\ga>1$.
\medskip

{\bf Condition C:}

More generally, let $\psi (t)$ be  strictly monotone increasing for $t $   large, with 
$\lim_{t\to\ff}\psi (t)=\ff$.  Assume that
\begin{equation}
\int_{0}^{1}z\,\psi (1/z)\,d\eta(z)<\ff.\label{3.80f}
\end{equation}
It follows as before that
\begin{equation}
|g(x)-g(y)|\leq  \frac{1}{c \psi (1/|x-y|)}.  \label{3.8lf}
\end{equation}
It follows that $d(x,y)\leq  \frac{1}{c \psi^{1/2} (1/|x-y|)}$, and 
then in the notation of \cite[Theorem 1.2]{MRsuf}
\bea
B_{d}(t,u)& =&\{x\,|\, \, d(x,t)\leq u\}\supseteq \{x\,|\, \frac{1}{c \psi (1/|x-t|)}\leq u^{2}\}\label{3.82f}\\
&=&   \{x\,|\, \psi (1/|x-t|)\geq \frac{1}{cu^{2}}\}\nn\\
&=&   \{x\,|\,  1/|x-t|\geq \psi^{-1}( 1/cu^{2})\}=  \{x\,|\, |x-t|\leq \frac {1}{\psi^{-1}( 1/cu^{2})}\}, \nn
\eea
so that if $\mu$ is a multiple of Lebesgue measure we have
\begin{equation}
\mu (B_{d}(t,u))\geq       \frac {1}{\psi^{-1}( 1/cu^{2})} \label{}
\end{equation}  
and therefore, by (\ref{3.90a}), we will have continuity if 
\begin{equation}
\lim_{\ep\to 0}\int_{0}^{\ep}  \(\log \psi^{-1}( 1/cu^{2})\)^{1/2}\,du=0.\label{3.83f}
\end{equation}

We have been unable to find necessary and sufficient conditions on $\eta$ for the joint continuity of $\wt L^{u}_{t}$  on $[0,1]\times R_{+}^{1}$.

See Appendix A
for a different way to do the algebra of (\ref{3.64})

\subsection{$\mathbf{p=0}$, $\mathbf{q=0}$}

In this case $u^{1}(y,0)=0$, so we cannot normalize the local time at $0$ as in (\ref{3.51}). In \cite[III, Theorem 3.7]{Blu} it is shown that there is a local time $\wh L^{0}_{t}$ at $0$ such that
\be
 \wt E^{ y}\(\int_{0}^{\ff}  e^{-\la t}\, d_{t}\wh L^{0}_{t}\)\label{b3.86b}
=e^{-\sqrt{2\la}\, y} \frac{1}{\int_{0}^{\ff} \(1-e^{-\sqrt{2\la}\,x}\)\,d\nu (x)}.
\ee
In particular, taking $\la=1$, and then using (\ref{3.51a}) we have, see (\ref{3.51}), 
\begin{equation}
  \wt E^{ y}\(\int_{0}^{\ff}  e^{- t}\, d_{t}\wh L^{0}_{t}\)= \frac{e^{-\sqrt{2}\, y}}{\int_{0}^{\ff} \(1-e^{-\sqrt{2}\,x}\)\,d\nu (x)}=e^{-\sqrt{2}\, y},\label{b3.87}
\end{equation}
while $u^{1}(y,0)=0$.
In \cite[p. 67]{Blu}, when $p=0$, there is a different presentation of the local time at $0$,
which is similar to that in \cite{IM1}, and which also covers the case where $q\geq 0$. This is described in Appendix B, where we also obtain (\ref{b3.86b}) when $q=0$.

On the other hand, for $u>0$  we have 
 \bea
\wt E^{ y}\(\int_{0}^{\ff}  e^{- t}\, d_{t}\wt L^{u}_{t}\)=u^{1}(y,u) &=&    v^{1}(y,u)  +e^{-\sqrt{2}\,\, y}\,\,\frac{  \int_{0}^{\ff}v^{1}(z,u) \,d\eta(z)}{\int_{0}^{\ff}\(1-e^{-\sqrt{2}\,\, z}\)\,d\eta(z)}\nn\\
&=&   v^{1}(y,u)  +e^{-\sqrt{2 }\,\, y}\,\,  \int_{0}^{\ff}v^{1}(z,u) \,d\eta(z).\label{blu.42denst}
\eea
Set
\begin{equation}
\wh L^{u}_{t}=\frac{\wt L^{u}_{t}}{  \int_{0}^{\ff}v^{1}(z,u) \,d\eta(z)},\label{3.88}
\end{equation}
so that by (\ref{blu.42denst})  we have 
 \begin{equation} \label{3.89}
  \wt E^{ y}\(\int_{0}^{\ff}  e^{- t}\, d_{t}\wh L^{u}_{t}\)= \left\{
 \begin{array} {cc}
\frac{v^{1}(y,u) }{ \int_{0}^{\ff}v^{1}(z,u) \,d\eta(z)}+e^{-\sqrt{2}\,\, y},& \quad u>0  
 \\
\hspace{1in}e^{-\sqrt{2}\, y},& \quad u=0   
\end{array}  \right. .
\end{equation}
Note that this is the potential density of $\wt X$ with respect to $f(u)\,du$ where 
$f(u)= \int_{0}^{\ff}v^{1}(z,u) \,d\eta(z)$.
We would like to show that $\wh L^{u}_{t}$ is     jointly continuous on $[0,1]\times R_{+}^{1}$. We provide three examples.
\medskip

{\bf Example 1:} Assume that 
\begin{equation}
d\eta(z)=\frac{1}{z^{1+\bb}}\,dz,\label{3.90s}
\end{equation}
for some $0<\bb<1 $. This corresponds to the stable subordinator of index $\bb$. Then with $u<1$,
\bea
&&\int_{0}^{\ff}v^{1}(z,u) \,d\eta(z)\geq \int_{0}^{u}v^{1}(z,u) \,d\eta(z)\label{3.91}\\
&&=\int_{0}^{u}e^{-\sqrt{2}\,\,u}\(\frac{e^{\sqrt{2}\,\,z}-e^{-\sqrt{2}\,\,z}}{\sqrt{2}}\) \,d\eta(z)\nn\\
&&=\int_{0}^{u}\frac{e^{-\sqrt{2}\,\,u}}{\sqrt{2}}\frac{\(e^{\sqrt{2}\,\,z}-e^{-\sqrt{2}\,\,z}\)}{z^{1+\bb}}\,dz \nn\\
&&\sim  2  \int_{0}^{u}e^{-\sqrt{2}\,\,u}\frac{z}{z^{1+\bb}}\,dz \sim \frac{2}{1-\bb}\,\, u^{1-\bb},\nn
\eea
as $u\to 0$.
It then follows from  (\ref{3base})     that
\begin{equation}
\frac{v^{1}(y,u) }{ \int_{0}^{\ff}v^{1}(z,u) \,d\eta(z)}\leq cu^{\bb},\label{3.93}
\end{equation}
which shows that (\ref{3.89}) is continuous in  $u$. In (\ref{3.91}) we gave a lower bound which will be adequate for our purposes, but in (\ref{3.105e7}) we show that it is actually the right order of magnitude.

We now show that 
  $\wh L^{u}_{t}$ is    jointly continuous on $[0,1]\times R_{+}^{1}$ for our example (\ref{3.90s}). If $\Gamma (y,u)$ denotes the right hand side of (\ref{3.89}) we have
  \begin{equation}
\sup_{y,u \in [0,1]}\Ga (y,u)<\ff\label{3.89a}
\end{equation}
and
\begin{equation}
\inf_{y,u \in [0,1]}\Ga (y,u)>0.\label{3.89b}
\end{equation}

We will show that for some $\bar \beta<\beta$
\begin{equation}
d^{2}(y,u)=\Ga^{2} (y,y)+\Ga^{2} (u,u)-2\Ga (y,u)\Ga (u,y)\leq c|y-u|^{\bar \bb}.\label{3.89u}
\end{equation}
  
Consider first the variation in $u$. 
\bea
&&\Gamma (y,u)-\Gamma (y,u')=\frac{v^{1}(y,u) }{ \int_{0}^{\ff}v^{1}(z,u) \,d\eta(z)}-\frac{v^{1}(y,u') }{ \int_{0}^{\ff}v^{1}(z,u') \,d\eta(z)}\label{3.94}\\
&& =         \frac{v^{1}(y,u)\int_{0}^{\ff}v^{1}(z,u') \,d\eta(z)- v^{1}(y,u')\int_{0}^{\ff}v^{1}(z,u) \,d\eta(z)}{ \int_{0}^{\ff}v^{1}(z,u) \,d\eta(z) \int_{0}^{\ff}v^{1}(z,u') \,d\eta(z)},       \nn
\eea
and write 
\begin{eqnarray}
&&v^{1}(y,u)\int_{0}^{\ff}v^{1}(z,u') \,d\eta(z)- v^{1}(y,u')\int_{0}^{\ff}v^{1}(z,u) \,d\eta(z)
\label{3.95}
\\
&&=\(v^{1}(y,u)-v^{1}(y,u')\)\int_{0}^{\ff}v^{1}(z,u') \,d\eta(z)\nn
\nonumber\\
&&+v^{1}(y,u')\(\int_{0}^{\ff}v^{1}(z,u') \,d\eta(z)-\int_{0}^{\ff}v^{1}(z,u) \,d\eta(z)\).\nn
\nonumber
\end{eqnarray}

We have that for any $\bb'>\bb$,  
\begin{equation}
\int_{0}^{1}z^{\bb'}\,d\eta(z)<\ff,\label{3.97}
\end{equation}
so that 
\begin{equation}
\int_{0}^{\ff}z^{\bb'}\wedge 1\,d\eta(z)=C_{\bb'}<\ff.\label{3.98}
\end{equation}
  Then by (\ref{3.501b}) 
\begin{equation}
|    \int_{0}^{\ff}\(v^1 (z,u)-v^1 (z,u')\)\,d\eta(z)      |\leq  4C_{\bb'}|u-u'|^{1-\bb'}.
\label{3.99}
\end{equation}

Then, using  (\ref{3.501a}) and (\ref{3base})
 we see that
\begin{eqnarray}
&&\Big | \frac{v^{1}(y,u) }{ \int_{0}^{\ff}v^{1}(z,u) \,d\eta(z)}-\frac{v^{1}(y,u') }{ \int_{0}^{\ff}v^{1}(z,u') \,d\eta(z)}\Big |
\label{3.100}
\\
&&\leq c\frac{|u-u'|}{ \int_{0}^{\ff}v^{1}(z,u) \,d\eta(z)}+c \frac{u'\, |u-u'|^{1-\bb'}}{ \int_{0}^{\ff}v^{1}(z,u') \,d\eta(z) \int_{0}^{\ff}v^{1}(z,u) \,d\eta(z)}
\nonumber\\
&&\leq c\frac{|u-u'|}{u^{1-\bb}}+c \frac{u'\, |u-u'|^{1-\bb'}}{ u^{1-\bb}u'^{1-\bb}}=
c\frac{|u-u'|}{u^{1-\bb}}+c \frac{u'^{\bb}\, |u-u'|^{1-\bb'}}{ u^{1-\bb}}.
\nonumber
\end{eqnarray}

If $u\geq u'$ this is bounded by      $c\frac{|u-u'|}{u^{1-\bb}}+c \frac{ |u-u'|^{1-\bb'}}{ u^{1-2\bb}}$, hence if $\bb\leq 1/2$, using     $u\geq |u-u'|$   we obtain the bound       $c|u-u'|^{\bb}+c|u-u'|^{2\bb-\bb'} $. If $\bb>1/2$ then $\bb=1-\bb +\al$ for some $\al>0$
and we have $u'^{\bb}=u'^{1-\bb}u'^{\al}\leq u^{1-\bb}$ using the fact that $u\geq u'$ and $u'\leq 1. $ In this case we have the bound $c|u-u'|^{\bb}+c|u-u'|^{1-\bb'} $. Thus for any $0<\bb<1$ we have the bound $c|u-u'|^{(2\bb)\wedge 1-\bb'} $.

In case $u'\geq u$ we can decompose 
(\ref{3.95}) differently:
\begin{eqnarray}
&&v^{1}(y,u)\int_{0}^{\ff}v^{1}(z,u') \,d\eta(z)- v^{1}(y,u')\int_{0}^{\ff}v^{1}(z,u) \,d\eta(z)
\label{3.101}
\\
&&=-\(v^{1}(y,u')-v^{1}(y,u)\)\int_{0}^{\ff}v^{1}(z,u) \,d\eta(z)\nn
\nonumber\\
&&+v^{1}(y,u)\(\int_{0}^{\ff}v^{1}(z,u') \,d\eta(z)-\int_{0}^{\ff}v^{1}(z,u) \,d\eta(z)\).\nn
\nonumber
\end{eqnarray}
Thus $|\Gamma (y,u)-\Gamma (y,u')|$ is bounded by $c|u-u'|^{(2\bb)\wedge 1-\bb'} $, so that
\begin{equation}
|\Gamma (y,y)-\Gamma (y,u)|\leq c|y-u|^{(2\bb)\wedge 1-\bb'}.\label{3.102s}
\end{equation}

Now consider the variation in $y$.
\begin{eqnarray}
&&|\Gamma (y,u)-\Gamma (y',u)|=\frac{|v^{1}(y,u)-v^{1}(y',u) |}{ \int_{0}^{\ff}v^{1}(z,u) \,d\eta(z)}+|e^{-\sqrt{2}\,\, y}-e^{-\sqrt{2}\,\, y'}|
\label{3.102}
\\
&&\leq c\frac{|v^{1}(y,u)-v^{1}(y',u) |}{u^{1-\bb}}+c|y-y'|.
\nonumber
\end{eqnarray}
But
\bea
&& |v^{1}(y,u)-v^{1}(y',u) |=|v^{1}(y,u)-v^{1}(y',u) |^{\bb}      |v^{1}(y,u)-v^{1}(y',u) |^{1-\bb}\label{3.103}\\
&&\leq c |y-y'|^{\bb}   |v^{1}(y,u)+v^{1}(y',u) |^{1-\bb}\leq   c |y-y'|^{\bb}  u^{1-\bb}.\nn
\eea
It follows that $|\Gamma (y,u)-\Gamma (y',u)|\leq c |y-y'|^{\bb}$ so that
\begin{equation}
|\Gamma (y,y)-\Gamma (u,y)|\leq c|y-u|^{\bb}.\label{3.102t}
\end{equation}
Together with (\ref{3.102}), (\ref{3.64}) and (\ref{3.89a}) this gives (\ref{3.89u}).
It follows from \cite[Theorem 1.2]{MRsuf} as before   that $\wh L^{u}_{t}$ is  jointly continuous on $[0,1]\times R_{+}^{1}$.
\medskip

We have not been able to find an exact uniform modulus of continuity  for  $\wh L^{u}_{t}$  on $[0,1]\times R_{+}^{1}$, or even an exact law of the iterated logarithm at $0$.
However it follows from \cite[Theorem 1.3]{MRsam} and its adoption to local times as in \cite{MRsuf} that for some $\bar\beta<\bb$
 \begin{equation}
\limsup_{h\to 0}\sup_{\stackrel{ x,y\in [0,1] }{|x-y|\le h }}\frac{|\wh L_{t}^{x}-\wh L_{t}^{y}|}{h^{\bar\beta/2} }=0, \quad a.e.\quad t, \quad a.s.\label{3.105}
\end{equation} 

  In particular, if $\bb\leq 1/2$ then (\ref{3.105}) holds for any $\bar\beta<\bb$, while if $\bb> 1/2$ then (\ref{3.105}) holds for any $\bar\beta<1-\bb$.

We can also obtain Holder continuity: for all  $x,y\in [0,1]$
 \begin{equation}
\Big  | \wh L_{t}^{x}- \wh L_{t}^{y}\,\,\Big |\leq C'(\om )\,\,|x-y|^{\bar\beta/2}, \quad \quad \forall t,  \label{3.105h2}
\end{equation} 
for some random variable $C'(\om )<\ff$, a.s. with $\bar\beta$ as above.

Proof of (\ref{3.105h2}): We first show that
 \begin{equation}
\Big  |\sqrt{\wh L_{t}^{x}}-\sqrt{\wh L_{t}^{y}}\,\,\Big |\leq C(\om )|x-y|^{\bar\beta/2}, \quad \quad \forall t,  \label{3.105h}
\end{equation} 
for some random variable $C(\om )<\ff$, a.s. with $\bar\beta$ as above. Using $a-b=(\sqrt{a}-\sqrt{b})(\sqrt{a}+\sqrt{b})$ for $a,b\geq 0$ and the fact that $\sup_{u\in [0,1]}\sqrt{\wh L_{t}^{u}}<\ff$, a.s. this implies (\ref{3.105h2}).

Let $\theta_{u}$ be the $1/2$-permanental process whose kernel is the right hand side of (\ref{3.89}). By \cite[(35) and Lemma 3.1]{MRsuf} we have that for any $x,y\in [0,1]$ and any $n$
\begin{equation}
E\(\(\theta^{1/2}_{x}-\theta^{1/2}_{y}\)^{2n}\)\leq \(cd^{2}(x,y)\)^{n}\leq  \(c|x-y|^{\bar\beta}\)^{n}.\label{3.105h3}
\end{equation}
It then follows from \cite[Theorem 14.1.1]{book} that for any $\al<\frac{n\bar\beta+1 }{2n}$
\begin{equation}
\Big | \theta^{1/2}_{x}-\theta^{1/2}_{y}\,\,\Big |\leq C(\om )|x-y|^{\al},\label{3.105h}
\end{equation}
for some $C(\om )<\ff$ a.s. By the Isomorphism Theorem and the joint continuity of the local times we obtain (\ref{3.105h}) first with    $\alpha$ instead of $\bar\beta/2$    and 
then the version we want by taking $n$ sufficiently large.

\medskip

{\bf Example 2:} Assume now that
  \begin{equation}
d\eta(z)=\frac{S(z)}{z^{1+\bb}}\,dz,\label{3.106}
\end{equation}
for some $0<\bb<1$,  where $S(z)$ is slowly varying at $0$ and bounded away from $0$.
In this case with $u<1$,
\bea
&&\int_{0}^{\ff}v^{1}(z,u) \,d\eta(z)\geq \int_{0}^{u}v^{1}(z,u) \,d\eta(z)\label{3.108}\\
&&=\int_{0}^{u}e^{-\sqrt{2}\,\,u}\(\frac{e^{\sqrt{2}\,\,z}-e^{-\sqrt{2}\,\,z}}{\sqrt{2}}\) \,d\eta(z)\nn\\
&&=\int_{0}^{u}\frac{e^{-\sqrt{2}\,\,u}}{\sqrt{2}}\frac{\(e^{\sqrt{2}\,\,z}-e^{-\sqrt{2}\,\,z}\)}{z^{1+\bb}}S(z)\,dz \nn\\
&&\sim  2  \int_{0}^{u}e^{-\sqrt{2}\,\,u}\frac{zS(z)}{z^{1+\bb}}\,dz \sim \frac{2}{1-\bb}\,\, S(u) u^{1-\bb},\nn
\eea
as $u\to 0$.

Everything should go through as before, but we also want to consider the case where instead of $z^{\bb'}$ with $\bb'>\bb$,  we use $z^{\bb}\bar S(z)$ where  $\bar S(z)$ is a a function which is slowly varying at $0$, bounded away from $0$ and such that
\begin{equation}
\int_{0}^{u}z^{\bb}\bar S(z)\,d\eta(z)=\int_{0}^{u}\frac{S(z)\bar S(z)}{z}\,dz <\ff,\label{3.109}
\end{equation}
for $u$ small
so that 
\begin{equation}
\int_{0}^{\ff}z^{\bb}\bar S(z)\wedge 1\,d\eta(z)=C<\ff.\label{3.110}
\end{equation}
Let $A(z)=z^{\bb}\bar S(z)$ and $B(z)=\frac{z^{1-\bb}}{\bar S(z)} $.
Then instead of (\ref{3.99}) we write 
\begin{eqnarray}
&&  |    \int_{0}^{\ff}\(v^1 (z,u)-v^1 (z,u')\)\,d\eta(z)      |
\label{3.111}
\\  
&& \leq  \int_{0}^{\ff}|v^1 (z,u)-v^1 (z,u')|\,d\eta(z)
\nonumber\\  
&& \leq  \int_{0}^{\ff} B\(|v^1 (z,u)-v^1 (z,u')|\)  A\(|v^1 (z,u)-v^1 (z,u')|\)\,d\eta(z)
\nonumber\\  
&& \leq  4 \frac{|u-u'|^{1-\bb}}{\bar S(|u-u'|)}     \int_{0}^{\ff}z^{\bb}\bar S(z)\wedge 1\,d\eta(z)=C\frac{|u-u'|^{1-\bb}}{\bar S(|u-u'|)} .
\nonumber
\end{eqnarray}
As in (\ref{3.100}) we have that 
\begin{eqnarray}
&&\Big | \frac{v^{1}(y,u) }{ \int_{0}^{\ff}v^{1}(z,u) \,d\eta(z)}-\frac{v^{1}(y,u') }{ \int_{0}^{\ff}v^{1}(z,u') \,d\eta(z)}\Big |
\label{3.112}
\\
&&\leq c\frac{|u-u'|}{ \int_{0}^{\ff}v^{1}(z,u) \,d\eta(z)}+c \frac{u'\, |u-u'|^{1-\bb}/\bar S(|u-u'|)}{ \int_{0}^{\ff}v^{1}(z,u') \,d\eta(z) \int_{0}^{\ff}v^{1}(z,u) \,d\eta(z)}
\nonumber\\
&&\leq c\frac{|u-u'|}{S(u)u^{1-\bb}}+c \frac{u'\, |u-u'|^{1-\bb}/\bar S(|u-u'|)}{S(u) u^{1-\bb}S(u')u'^{1-\bb}}\nonumber\\
&&=c\frac{|u-u'|}{S(u)u^{1-\bb}}+c \frac{u'^{\bb}\, |u-u'|^{1-\bb}/\bar S(|u-u'|)}{S(u')S(u) u^{1-\bb}}.
\nonumber
\end{eqnarray}
If $u\geq u'$ and    $\bb\leq 1/2$   we have $u\geq  |u-u'|$ and also $\frac{u'^{\bb}}{S(u')}\leq   \frac{u^{\bb}}{S(u)}$ and find that the last display is bounded by
\begin{equation}
c\frac{|u-u'|^{\bb}}{S(|u-u'|)}+c \frac{|u-u'|^{\bb}}{S(|u-u'|)S(|u-u'|)\bar S(|u-u'|)},\label{}
\end{equation}
with an analogous result for $\bb>1/2$, and similar results for $u'\geq u$. By (\ref{3.109}) we must have $S(|u-u'|)\bar S(|u-u'|)\to 0$ so  we will have 
\begin{equation}
d^{2}(y,u)\leq   c \frac{|u-u'|^{\bb \wedge (1-\bb)}}{S^{2}(|u-u'|)\bar S(|u-u'|)}. \label{}
\end{equation}

{\bf Example 3:} We now consider the case where $\bb=1$. 
 Let 
   \begin{equation}
 d\eta(z)=\frac{|\log^{2m-1} (z)|}{z^{2}\(1+\log^{2m} (z)\)^{2}},\label{3.141}
 \end{equation}
 where $m$ is an integer $\geq 2$.
 We first verify that $\eta$ is a L\'evy measure. For this it suffices to show that
   \begin{equation}
\int_{0}^{\ff} (z\wedge 1) d\eta(z)=\int_{0}^{\ff} (z\wedge 1) \frac{|\log^{2m-1} (z)|}{z^{2}\(1+\log^{2m} (z)\)^{2}}\,dz<\ff.\label{3.141p}
\end{equation}
It is clear that the integral from $1$ to $\ff$ is finite. 
Note that
 \begin{equation}
\frac{d}{dz}\(\frac{1}{1+\log^{2m} (z)}\)=2m\frac{-\log^{2m-1} (z)}{z\(1+\log^{2m} (z)\)^{2}}.\label{3.140d}
 \end{equation} 
Since  $\log^{2m-1} (z)<0$ from $0$ to $1$, using  (\ref{3.140d})
\bea
&&\int_{0}^{1} (z\wedge 1) \frac{|\log^{2m-1} (z)|}{z^{2}\(1+\log^{2m} (z)\)^{2}}\,dz=
-\int_{0}^{1} \frac{\log^{2m-1} (z)}{z\(1+\log^{2m} (z)\)^{2}}\,dz\label{}\\
&&=\frac{1}{2m}\frac{1}{(1+\log^{2m} (z))}\Bigg|_{0}^{1}=\frac{1}{2m}.
\eea
Using  (\ref{3.140d}) again have for $u<1$
\bea
&&\int_{0}^{u}  \frac{|\log^{2m-1} (z)|}{z\(1+\log^{2m} (z)\)^{2}}=  \int_{0}^{u}  \frac{-\log^{2m-1} (z)}{z\(1+\log^{2m} (z)\)^{2}}   \label{3.142}\\
&&=\frac{1}{2m}  \frac{1}{1+\log^{2m}(z)}\Bigg |_{0}^{u}    =\frac{1}{2m}  \frac{1}{1+\log^{2m}(u)}.\nn
\eea
Hence, see (\ref{3.91}), 
\begin{equation}
\int_{0}^{\ff}v^{1}(z,u) \,d\eta(z)\geq \int_{0}^{u}v^{1}(z,u) \,d\eta(z)\geq  \frac{1}{m} \,\,  \log^{-2m}(u),\label{3.143}
\ee
as $u\to 0$.
It then follows from  (\ref{3base})     that
\begin{equation}
\frac{v^{1}(y,u) }{ \int_{0}^{\ff}v^{1}(z,u) \,d\eta(z)}\leq cu\log^{2m}(u),\label{3.144}
\end{equation}
as $u\to 0$, which shows that (\ref{3.89}) is continuous in  $u$.

We  also note that 
\begin{eqnarray}
&&\int_{0}^{1}z\,\log ^{3}(1/z)\,d\eta(z)= \int_{0}^{1} \frac{\log^{2m+2} (z)}{z\(1+\log^{2m} (z)\)^{2}}\,dz  
\label{3.145}
\\
&&\leq  \int_{0}^{1}    \frac{\log^{2} (z)}{z\(1+\log^{2m} (z)\)}\,dz<\ff,
\nonumber
\end{eqnarray}
since with $z=e^{y}$, $dz=e^{y}\,dy$ and 
\begin{equation}
\int_{0}^{1}    \frac{\log^{2} (z)}{z\(1+\log^{2m} (z)\)}\,dz=\int_{-\ff}^{0}    \frac{y^{2}  }{\(1+y^{2m}\)}\,dy<\ff,\label{3.146}
\end{equation}
for $m\geq 2$. Compare (\ref{3.51aal}),
so that with $h(t)=t\,\log_{\ast} ^{3}(1/t)$, where $\log_{\ast}(x)=\log(x)$ if $x\geq 1$ and 
$\log_{\ast}(x)=1$ if $x<1 $, 
\begin{equation}
\int_{0}^{\ff}h\(z\)\wedge 1\,d\eta(z)=c_{h}<\ff.\label{3.501ls}
\end{equation}
It then follows as before that  
\begin{equation}
|    \int_{0}^{\ff}\(v^1 (z,u)-v^1 (z,u')\)\,d\eta(z)      |\leq c \log_{\ast}  ^{-3}(1/|u-u'|).  \label{3.57ls}
\end{equation}

We now show that 
  $\wh L^{u}_{t}$ is    jointly continuous on $[0,1]\times R_{+}^{1}$ for our example (\ref{3.141}). If $\Gamma (y,u)$ denotes the right hand side of (\ref{3.89}) we have
  \begin{equation}
\sup_{y,u \in [0,1]}\Ga (y,u)<\ff\label{3.147}
\end{equation}
and
\begin{equation}
\inf_{y,u \in [0,1]}\Ga (y,u)>0.\label{3.148}
\end{equation}

We will show that  
\begin{equation}
d^{2}(y,u)=\Ga^{2} (y,y)+\Ga^{2} (u,u)-2\Ga (y,u)\Ga (u,y)\leq c \log_{\ast}^{-3}(1/|u-y|).\label{3.149}
\end{equation}
It will then follow as in example (\ref{3.51aal}) that   $\wh L^{u}_{t}$ is    jointly continuous on $[0,1]\times R_{+}^{1}$.
  
Consider first the variation in $u$. 
\bea
&&\Gamma (y,u)-\Gamma (y,u')=\frac{v^{1}(y,u) }{ \int_{0}^{\ff}v^{1}(z,u) \,d\eta(z)}-\frac{v^{1}(y,u') }{ \int_{0}^{\ff}v^{1}(z,u') \,d\eta(z)}\label{3.150}\\
&& =         \frac{v^{1}(y,u)\int_{0}^{\ff}v^{1}(z,u') \,d\eta(z)- v^{1}(y,u')\int_{0}^{\ff}v^{1}(z,u) \,d\eta(z)}{ \int_{0}^{\ff}v^{1}(z,u) \,d\eta(z) \int_{0}^{\ff}v^{1}(z,u') \,d\eta(z)},       \nn
\eea
and write 
\begin{eqnarray}
&&v^{1}(y,u)\int_{0}^{\ff}v^{1}(z,u') \,d\eta(z)- v^{1}(y,u')\int_{0}^{\ff}v^{1}(z,u) \,d\eta(z)
\label{3.151}
\\
&&=\(v^{1}(y,u)-v^{1}(y,u')\)\int_{0}^{\ff}v^{1}(z,u') \,d\eta(z)\nn
\nonumber\\
&&+v^{1}(y,u')\(\int_{0}^{\ff}v^{1}(z,u') \,d\eta(z)-\int_{0}^{\ff}v^{1}(z,u) \,d\eta(z)\).\nn
\nonumber
\end{eqnarray}

Putting this together we have that 
\begin{eqnarray}
&&\Big | \frac{v^{1}(y,u) }{ \int_{0}^{\ff}v^{1}(z,u) \,d\eta(z)}-\frac{v^{1}(y,u') }{ \int_{0}^{\ff}v^{1}(z,u') \,d\eta(z)}\Big |
\label{3.155}
\\
&&\leq c\frac{|u-u'|}{ \int_{0}^{\ff}v^{1}(z,u) \,d\eta(z)}+c \frac{u'\,  \log_{\ast}  ^{-3}(1/|u-u'|)}{ \int_{0}^{\ff}v^{1}(z,u') \,d\eta(z) \int_{0}^{\ff}v^{1}(z,u) \,d\eta(z)}
\nonumber\\
&&\leq c|u-u'|\log^{2m}(1/u)+c  u'\log^{2m}(1/u')\log^{2m}(1/u)\,  \log_{\ast}  ^{-3}(1/|u-u'|)\nonumber
\end{eqnarray}

If $u\geq u'$,  we have that  $u'\log^{2m}(1/u')\log^{2m}(1/u)\leq  u'\log^{4m}(1/u') \leq c $ and        using     $u\geq |u-u'|$   we obtain the bound $|u-u'|\log^{2m}(1/u)\leq |u-u'|\log^{2m}(1/|u-u'|)   $. Thus when     $u\geq u'$,     (\ref{3.155}) is bounded by $c\log_{\ast}  ^{-3}(1/|u-u'|)$, and as before when $u'\geq u$ we can decompose 
(\ref{3.151}) differently and obtain a similar bound.

Thus $|\Gamma (y,u)-\Gamma (y,u')|$ is bounded by $c\log_{\ast}  ^{-3}(1/|u-u'|)$, so that
\begin{equation}
|\Gamma (y,y)-\Gamma (y,u)|\leq c\log_{\ast}  ^{-3}(1/|y-u|).\label{3.157}
\end{equation}

Now consider the variation in $y$.
\begin{eqnarray}
&&|\Gamma (y,u)-\Gamma (y',u)|=\frac{|v^{1}(y,u)-v^{1}(y',u) |}{ \int_{0}^{\ff}v^{1}(z,u) \,d\eta(z)}+|e^{-\sqrt{2}\,\, y}-e^{-\sqrt{2}\,\, y'}|
\nn
\\
&&\leq c |v^{1}(y,u)-v^{1}(y',u) |\log^{2m}(1/u) +c|y-y'|.
\label{3.158}
\end{eqnarray}

But
\bea
&& |v^{1}(y,u)-v^{1}(y',u) |=|v^{1}(y,u)-v^{1}(y',u) |^{1/2}      |v^{1}(y,u)-v^{1}(y',u) |^{1/2}\nn\\
&&\leq c |y-y'|^{1/2}   |v^{1}(y,u)+v^{1}(y',u) |^{1/2}\leq   c |y-y'|^{1/2}  u^{1/2}.\label{3.159}
\eea
It follows that $|\Gamma (y,u)-\Gamma (y',u)|\leq c |y-y'|^{1/2}$ so that
\begin{equation}
|\Gamma (y,y)-\Gamma (u,y)|\leq c|y-u|^{1/2}.\label{3.160}
\end{equation}
As before this gives (\ref{3.149}). As mentioned,  this completes the proof that  $\wh L^{u}_{t}$ is    jointly continuous on $[0,1]\times R_{+}^{1}$.
\medskip

In (\ref{3.143}) we gave a lower bound for $\int_{0}^{\ff}v^{1}(z,u) \,d\eta(z)$   as $u\to 0$. In Appendix D we show that it is the right order of magnitude.


\begin{appendices}

\section{A different way to do the algebra of (\ref{3.64})}



We write out
\begin{eqnarray}
&&\Ga^{2} (x,x)+\Ga^{2} (y,y)-2\Ga (x,y)\Ga (y,x)
\label{3.60}
\\
&&=(v^1)^{2}(x,x)+2v^1(x,x)f(x)g(x) +f^{2}(x)g^{2}(x)
\nonumber\\
&&+(v^1)^{2}(y,y)+2v^1(y,y)f(y)g(y) +f^{2}(y)g^{2}(y)
\nonumber\\
&&-2    \( v^1(x,y)v^1(y,x)+v^1(x,y)f(y)g(x) +v^1(y,x)f(x)g(y)+f(x)g(x)f(y)g(y)\)
\nonumber\\
&&=(v^1)^{2}(x,x)+(v^1)^{2}(y,y) -2     v^1(x,y)v^1(y,x)
\nonumber\\
&&+2\( v^1(x,x)f(x)g(x)+v^1(y,y)f(y)g(y)-  v^1(x,y)f(y)g(x) -v^1(y,x)f(x)g(y)      \)
\nonumber\\
&&+    \( f(x)g(x)-f(y)g(y)\)^{2}=:A+B+C.
\nonumber
\end{eqnarray}
Using the symmetry of $v^1(x,y)$, and the fact that $v^1(x,y)\leq v^1(x,x)\wedge v^1(y,y)$
\begin{eqnarray}
&&A=\((v^1)^{2}(x,x)-(v^1)^{2}(x,y)\) +\((v^1)^{2}(y,y)-(v^1)^{2}(x,y)\)
\label{3.61}
\\
&&=\(v^1(x,x)-v^1(x,y)\)\(v^1(x,x)+v^1(x,y)\)\nn\\
&&+\(v^1(y,y)-v^1(x,y)\)\(v^1(y,y)+v^1(x,y)\).
\nonumber\\
&&\leq \(v^1(x,x)+v^1(y,y)-2v^1(x,y)\)2\sup_{x\in S} v^1(x,x),
\nonumber
\end{eqnarray}

\begin{eqnarray}
&&B/2=v^1(x,x)f(x)g(x)+v^1(y,y)f(y)g(y)-  v^1(x,y)\(f(y)g(x) +f(x)g(y)\) 
\nn
\\
&&=\(v^1(x,x)+v^1(y,y)-2v^1(x,y)\)f(x)g(x)+v^1(y,y)(f(y)g(y)-f(x)g(x))
\nonumber\\
&& +v^1(x,y)\(2f(x)g(x)-f(y)g(x) -f(x)g(y)\),
\label{3.62}
\end{eqnarray}
and
\begin{eqnarray}
&&C=\( f(x)g(x)-f(y)g(y)\)^{2}
\label{3.63}
\\
&&=\( (f(x)g(x)-f(x)g(y))-(f(y)g(y))-f(x)g(y))\)^{2}
\nonumber\\
&&\leq 2\( f(x)g(x)-f(x)g(y)\)^{2}+2\((f(y)g(y)-f(x)g(y))\)^{2}
\nonumber\\
&&=2f^{2}(x)\( g(x)-g(y)\)^{2}+2g^{2}(y)\((f(y)-f(x))\)^{2}.
\nonumber
\end{eqnarray}


\section{The normalization in \cite{Blu}\\ for $\mathbf{p=0,\,q\geq 0}$}


In \cite[p. 67]{Blu}, when $p=0$, the local time at $0$, which we will denote by  $\wh L^{0}_{t}$, can be written as 
\begin{equation}
\wh L^{0}_{t}=T^{-1}\(\ell^{0}_{t}\),\label{3.80}
\end{equation}
where $\ell^{0}_{t}$ is the local time at $0$ of Brownian motion $W$ and 
\begin{equation}
T^{-1}(s)=\inf\{r\geq 0\,|\, T(r)>s\},\label{3.81}
\end{equation}
where $T(r),\, r\geq 0$ is the subordinator with
\begin{equation}
 E_{T}\(e^{-a \, T(r) }\)=e^{-r\(\frac{q}{\sqrt{2}}a+\int_{0}^{\ff} \(1-e^{-ax}\)\,d\nu (x)\)}, \label{3.82}
\end{equation}
and $T$ is independent of the Brownian motion $W$. We also recall that if 
\begin{equation}
\tau (s)=\inf\{t\geq 0\,|\, \ell^{0}_{t}>s\}\label{3.83}
\end{equation}
is the Brownian inverse local time, then 
\begin{equation}
E^{0}_{W}\(e^{-\la \tau (s)}\)=e^{-\sqrt{2\la}\, s}.\label{3.84}
\end{equation}

Using \cite[Theorem A4.3]{S} for change of variables,  it follows from all this that 
\begin{eqnarray}
&&E_{T}E^{0}_{W}\(\int_{0}^{\ff}e^{-\la t}dT^{-1}\(\ell^{0}_{t}\)\)=E_{T}E^{0}_{W}\(\int_{0}^{\ff}e^{-\la \tau (s)}dT^{-1}\(s\)\)
\label{3.85}
\\
&&=E_{T}\(\int_{0}^{\ff}E^{0}_{W}\(e^{-\la \tau (s)}\)dT^{-1}\(s\)\)=E_{T}\(\int_{0}^{\ff} e^{-\sqrt{2\la}\,\, s}\,dT^{-1}\(s\)\)
\nonumber\\
&&=E_{T}\(\int_{0}^{\ff} e^{-\sqrt{2\la}\,\, T\(r\)}\,dr\)=\int_{0}^{\ff} E_{T}\(e^{-\sqrt{2\la} \, T\(r\) }\)\,\,dr
\nonumber\\
&&=\int_{0}^{\ff}  e^{-r\(q\sqrt{\la}+\int _{0}^{\ff} \(1-e^{-\sqrt{2\la}\,x}\)\,d\nu (x)\)}\,\,dr=\frac{1}{q\sqrt{\la}+\int_{0}^{\ff} \(1-e^{-\sqrt{2\la}\,x}\)\,d\nu (x)}.
\nonumber
\end{eqnarray}

As explained in \cite[bottom of p. 64]{Blu}, if we start our Feller Brownian motion $\wt   X$ at some $y>0$, it behaves like Brownian motion until it's first return to $0$.
From this we see that for any $y\geq 0$
\begin{eqnarray}
&& \wt E^{ y}\(\int_{0}^{\ff}  e^{-\la t}\, d_{t}\wh L^{0}_{t}\)\label{3.86}
\\
&&=E_{T}E^{y}_{W}\(\int_{0}^{\ff}e^{-\la t}dT^{-1}\(\ell^{0}_{t}\)\)=e^{-\sqrt{2\la}\, y}E_{T}E^{0}_{W}\(\int_{0}^{\ff}e^{-\la t}dT^{-1}\(\ell^{0}_{t}\)\)
\nn
\\
&&=e^{-\sqrt{2\la}\, y} \frac{1}{q\sqrt{\la}+\int_{0}^{\ff} \(1-e^{-\sqrt{2\la}\,x}\)\,d\nu (x)}.
\nonumber
\end{eqnarray}

Comparing this with (\ref{blu.42dens0}) we see that when $q>0$
\begin{equation}
\wt L^{0}_{t}=q\sqrt{2}\,\,\wh L^{0}_{t}.\label{3.86a}
\end{equation}
When $q=0$ we have 
\be
 \wt E^{ y}\(\int_{0}^{\ff}  e^{-\la t}\, d_{t}\wh L^{0}_{t}\)\label{3.86b}
=e^{-\sqrt{2\la}\, y} \frac{1}{\int_{0}^{\ff} \(1-e^{-\sqrt{2\la}\,x}\)\,d\nu (x)},
\ee
while by (\ref{blu.42dens0}) we have $u^{\la}(y,0)=0$.
In particular, taking $\la=1$, and then using (\ref{3.51a}) we have, see (\ref{3.51}), 
\begin{equation}
  \wt E^{ y}\(\int_{0}^{\ff}  e^{- t}\, d_{t}\wh L^{0}_{t}\)= \frac{e^{-\sqrt{2}\, y}}{\int_{0}^{\ff} \(1-e^{-\sqrt{2}\,x}\)\,d\nu (x)}=e^{-\sqrt{2}\, y},\label{3.87}
\end{equation}
while $u^{1}(y,0)=0$.

\section{The order of magnitude of\\ $\mathbf{\int_{0}^{\ff}v^{1}(z,u) \,d\eta(z)}$ in (\ref{3.91})}

 
 In (\ref{3.91}) we gave a lower bound for $\int_{0}^{\ff}v^{1}(z,u) \,d\eta(z)$   as $u\to 0$. Here we show that it is the right order of magnitude.
\be
\int_{0}^{1}v^{1}(z,u) \,d\eta(z)= \int_{0}^{u}v^{1}(z,u) \,d\eta(z)+\int_{u}^{1}v^{1}(z,u) \,d\eta(z)\label{3.105e1}
\ee
and we saw that 
\begin{equation}
 \int_{0}^{u}v^{1}(z,u) \,d\eta(z)\sim \frac{2}{1-\bb}\,\, u^{1-\bb},\label{3.105e2}
\end{equation}
as $u\to 0$. We have 
\bea
&&\int_{u}^{1}v^{1}(z,u) \,d\eta(z)\label{3.105e3}\\
&&=\int_{u}^{1}e^{-\sqrt{2}\,\,z}\(\frac{e^{\sqrt{2}\,\,u}-e^{-\sqrt{2}\,\,u}}{\sqrt{2}}\) \,d\eta(z)= \(\frac{e^{\sqrt{2}\,\,u}-e^{-\sqrt{2}\,\,u}}{\sqrt{2}}\)\int_{u}^{1}e^{-\sqrt{2}\,\,z} \,\frac{1}{z^{1+\bb}}\,dz\nn\\
&&\sim    2u\int_{u}^{1}\(1-\sqrt{2}O(z)\) \frac{1}{z^{1+\bb}}\,dz=  2 u \int_{u}^{1} \(\frac{1}{z^{1+\bb}}-\sqrt{2}\frac{O(z)}{z^{1+\bb}}\)\,dz.\nn
\eea
We have 
\begin{equation}
2 u \int_{u}^{1}  \frac{1}{z^{1+\bb}}\,dz=2u\(-\frac{1}{\bb}\frac{1}{z^{\bb}}\Big |_{u}^{1}\)=-2u\frac{1}{\bb}+\frac{2}{\bb}u^{1-\bb}\label{3.105e4}
\end{equation}
while
\begin{equation}
2\sqrt{2} u \int_{u}^{1}  \frac{z}{z^{1+\bb}}\,dz=2\sqrt{2} u \int_{u}^{1}  \frac{1}{z^{\bb}}\,dz=2\sqrt{2} u\(\frac{1}{1-\bb} z^{1-\bb}\Big |_{u}^{1}\)= \frac{2\sqrt{2}}{1-\bb}\(u-u^{2-\bb}\)\label{3.105e5}
\end{equation}
as $u\to 0$. Of course
\begin{equation}
\int_{1}^{\ff}v^{1}(z,u) \,d\eta(z)=\(\frac{e^{\sqrt{2}\,\,u}-e^{-\sqrt{2}\,\,u}}{\sqrt{2}}\)\int_{1}^{\ff}e^{-\sqrt{2}\,\,z} \,\frac{1}{z^{1+\bb}}\,dz\sim cu\label{3.105e6}
\end{equation}
as $u\to 0$.

Thus we see that  
\begin{equation}
\int_{0}^{\ff}v^{1}(z,u) \,d\eta(z)\sim   \( \frac{2}{1-\bb}+   \frac{2}{\bb}\)u^{1-\bb}\quad \mbox{ as } u\to 0. \label{3.105e7}
\end{equation}  

\section{The order of magnitude of\\ $\mathbf{\int_{0}^{\ff}v^{1}(z,u) \,d\eta(z)}$ in (\ref{3.143})}

In (\ref{3.143}) we gave a lower bound for $\int_{0}^{\ff}v^{1}(z,u) \,d\eta(z)$   as $u\to 0$. We showed that
\begin{equation}
 \int_{0}^{u}v^{1}(z,u) \,d\eta(z)\sim \frac{1}{m}\,\, \log ^{-2m} (u), \quad \mbox{as $u\to 0$,}\label{3.105em2}
\end{equation}
 so that
\begin{equation}
 \int_{0}^{\ff}v^{1}(z,u) \,d\eta(z)\geq \frac{1}{m}\,\, \log ^{-2m} (u), \quad \mbox{as $u\to 0$.}\label{3.105em2}
\end{equation}
Here we show that this is the right order of magnitude. We have
\be
\int_{0}^{\ff}v^{1}(z,u) \,d\eta(z)= \int_{0}^{u}v^{1}(z,u) \,d\eta(z)+\int_{u}^{\ff}v^{1}(z,u) \,d\eta(z).\label{3.105em1}
\ee
 Hence it suffices to show that 
\begin{equation}
 \int_{u}^{\ff}v^{1}(z,u) \,d\eta(z)\leq c\,\, \log ^{-2m} (u), \quad \mbox{as $u\to 0$.}\label{bas.1}
\end{equation}

We have 
\bea
&&\int_{u}^{\ff}v^{1}(z,u) \,d\eta(z)\label{3.105em3}\\
&&=\int_{u}^{\ff}e^{-\sqrt{2}\,\,z}\(\frac{e^{\sqrt{2}\,\,u}-e^{-\sqrt{2}\,\,u}}{\sqrt{2}}\) \,d\eta(z)= \(\frac{e^{\sqrt{2}\,\,u}-e^{-\sqrt{2}\,\,u}}{\sqrt{2}}\)\int_{u}^{\ff}e^{-\sqrt{2}\,\,z} \,\,d\eta(z)\nn\\
&&\leq    2u\int_{u}^{\ff} \frac{|\log^{2m-1} (z)|}{z^{2}\(1+\log^{2m} (z)\)^{2}}\,dz, \quad \mbox{as $u\to 0$.}\nn
\eea
Clearly, $\int_{1}^{\ff} \frac{|\log^{2m-1} (z)|}{z^{2}\(1+\log^{2m} (z)\)^{2}}\,dz<\ff$ so that it suffices to show that
\begin{equation}
 \int_{u}^{1}v^{1}(z,u) \,d\eta(z)\leq c\,\, \log ^{-2m} (u), \quad \mbox{as $u\to 0$.}\label{bas.2}
\end{equation}

Using (\ref{3.140d}) which says that
 \begin{equation}
\frac{d}{dz}\(\frac{1}{1+\log^{2m} (z)}\)=2m\frac{-\log^{2m-1} (z)}{z\(1+\log^{2m} (z)\)^{2}},\label{3.140ds}
 \end{equation} 
we see that
\bea
&&2 u \int_{u}^{1}  \frac{|\log^{2m-1} (z)|}{z^{2}\(1+\log^{2m} (z)\)^{2}}\,dz=-2 u \int_{u}^{1}  \frac{\log^{2m-1} (z)}{z^{2}\(1+\log^{2m} (z)\)^{2}}\,dz\label{3.105em4}\\
&&=\frac{u}{m}  \int_{u}^{1}\frac{1}{z}\frac{d}{dz}\(\frac{1}{1+\log^{2m} (z)}\)\,dz    \nn\\
&&= \frac{u}{m} \,\,   \frac{1}{z\(1+\log^{2m} (z)\)} \Bigg |^{1}_{u}+\frac{u}{m}  \int_{u}^{1} \frac{1}{z^2\(1+\log^{2m} (z)\)}\,dz. 
  \nn
\eea
We have
\begin{equation}
\frac{u}{m} \,\,   \frac{1}{z\(1+\log^{2m} (z)\)} \Bigg |^{1}_{u}=\frac{u}{m}-\frac{1}{m} \frac{1}{\(1+\log^{2m} (u)\)}, \label{3.105em5}
\end{equation}
and write
\bea
&&\frac{u}{m}  \int_{u}^{1} \frac{1}{z^2\(1+\log^{2m} (z)\)}\,dz=\frac{u}{m} \int_{u}^{\log^{2m} (u)u} \frac{1}{z^2\(1+\log^{2m} (z)\)}\,dz\nn\\
&&\hspace{1.5 in}+\frac{u}{m}\int_{\log^{2m} (u)u}^{1} \frac{1}{z^2\(1+\log^{2m} (z)\)}\,dz.\nn
\eea
We have
\begin{equation}
\frac{u}{m}\int_{\log^{2m} (u)u}^{1} \frac{1}{z^2\(1+\log^{2m} (z)\)}\,dz\leq \frac{u}{m}\int_{\log^{2m} (u)u}^{1} \frac{1}{z^2}\,dz\leq \frac{1}{m} \frac{1}{\log^{2m} (u)}\label{3.105em16}
\end{equation}
while
\bea
&&\frac{u}{m}\int_{u}^{\log^{2m} (u)u} \frac{1}{z^2\(1+\log^{2m} (z)\)}\,dz\
\label{3.105em6}\\
&&\leq \frac{1}{\(1+\log^{2m} (\log^{2m} (u)u)\)}\frac{u}{m}  \int_{u}^{\log^{2m} (u)u}\frac{1}{z^2} \,dz\,\,\nn\\
&&\leq\frac{1}{m} \frac{1}{\(1+\log^{2m} (\log^{2m} (u)u)\)}.\nn
\eea

Combining the last three displays we see that
\begin{equation}
\frac{u}{m}  \int_{u}^{1} \frac{1}{z^2\(1+\log^{2m} (z)\)}\,dz\leq \frac{1}{m}  \frac{1}{\log^{2m} (u)}+ \frac{1}{m} \frac{1}{\(1+\log^{2m} (\log^{2m} (u)u)\)}.\label{3.105em18}
\end{equation}

Combined with (\ref{3.105em4}) and (\ref{3.105em5}) we see that  
\begin{equation}
\int_{u}^{1}v^{1}(z,u) \,d\eta(z)\leq \frac{1}{m} \frac{1}{\(1+\log^{2m} (\log^{2m} (u)u)\)}\asymp\,\, \log ^{-2m} (u)\quad \mbox{ as } u\to 0, \label{3.105em19}
\end{equation}   
which proves (\ref{bas.2}). 

\end{appendices}

\bibliographystyle{amsplain}

\bigskip
\noindent
\begin{tabular}{lll} &   P.J. Fitzsimmons  \\
&Department of
Mathematics  \\
& University of California,  \\
& San Diego,
La Jolla CA, USA  \\ 
& pfitzsim@ucsd.edu\\ 
& &\\
& & \\
& Jay Rosen\\
& Department of Mathematics\\
&  College of Staten Island, CUNY\\
& Staten Island, NY 10314, USA \\
& jrosen30@optimum.net
\end{tabular}

\end{document}